\documentclass[11pt]{article}

\usepackage{hyperref}
\usepackage{amsmath}
\usepackage{amssymb}
\usepackage{amsthm}
\usepackage{latexsym}
\usepackage{color}
\usepackage{graphicx}
\usepackage{appendix}
\usepackage{color} 
\usepackage{enumerate}

\usepackage[T1]{fontenc}
\usepackage[english]{babel}
\usepackage{geometry}
\DeclareSymbolFont{calletters}{OMS}{cmsy}{m}{n}
\DeclareSymbolFontAlphabet{\mathcal}{calletters}

\def\be{\begin{eqnarray}}
\def\ee{\end{eqnarray}}

\def\b*{\begin{eqnarray*}}
\def\e*{\end{eqnarray*}}

\newtheorem{Theorem}{Theorem}[part]
\newtheorem{Definition}{Definition}[part]
\newtheorem{Proposition}{Proposition}[part]

\newtheorem{Assumption}{Assumption}[part]
\newtheorem{Lemma}{Lemma}[part]
\newtheorem{Corollary}{Corollary}[part]
\newtheorem{Remark}{Remark}[part]
\newtheorem{Example}{Example}[part]

\makeatletter \@addtoreset{equation}{section}

\@addtoreset{Definition}{section}

\@addtoreset{Theorem}{section}

\@addtoreset{Proposition}{section}

\@addtoreset{Property}{section}

\@addtoreset{Assumption}{section}

\@addtoreset{Corollary}{section}

\@addtoreset{Lemma}{section}

\@addtoreset{Remark}{section}

\@addtoreset{Example}{section}

\newcommand{\No}[1]{\left\|#1\right\|}     
\newcommand{\abs}[1]{\left|#1\right|}     

\def \E{\mathbb{E}}

\def \L{\mathbb{L}}

\def \P{\mathbb{P}}
\def \Q{\mathbb{Q}}
\def \R{\mathbb{R}}

\def\Ac{{\cal A}}
\def\Bc{{\cal B}}

\def\Ec{{\cal E}}
\def\Fc{{\cal F}}

\def\Uh{{\hat U}}

\def\Zh{{\hat Z}}

\def\Ut{{\widetilde U}}

\def\Yt{{\widetilde Y}}
\def\Zt{{\widetilde Z}}

\def\esup{{\rm ess \, sup}}

\def\={\;=\;}
\def\.{\;.}

\def\eps{\varepsilon}

\def\reff#1{{\rm(\ref{#1})}}

\def\1{{\bf 1}}
\def \ep{\hbox{ }\hfill{ ${\cal t}$~\hspace{-5.1mm}~${\cal u}$   } }

\def \proof{{\noindent \bf Proof. }}
\def \ep{\hbox{ }\hfill$\Box$}

 \def\normeL2#1{\left\|{#1}\right\|_{L^2}}

 \title{Quadratic BSDEs with jumps: related non-linear expectations\footnote{Research supported by the Chair {\it Financial Risks} of the {\it Risk Foundation} sponsored by Soci\'et\'e G\'en\'erale, the Chair {\it Derivatives of the Future} sponsored by the {F\'ed\'eration Bancaire Fran\c{c}aise}, and the Chair {\it Finance and Sustainable Development} sponsored by EDF and Calyon.}
}
\author{Nabil {\sc Kazi-Tani}\footnote{CMAP, Ecole Polytechnique, Paris, nabil.kazitani@polytechnique.edu.} \and Dylan {\sc Possama\"{i}}\footnote{CEREMADE, Universit\'e Paris Dauphine, possamai@ceremade.dauphine.fr.}
      \and Chao {\sc Zhou}\footnote{Department of Mathematics, National University of Singapore, Singapore, matzc@nus.edu.sg. Part of this work was carried out while the author was working at CMAP, Ecole Polytechnique,  whose financial support is kindly acknowledged.} }          
 \date{\today}

\begin{document}

 \maketitle

\vspace{3mm}

 \begin{abstract}
\vspace{10mm}

In this article, we follow the study of quadratic backward SDEs with jumps,that is to say for which the generator has quadratic growth in the variables $(z,u)$, started in our accompanying paper \cite{kpz_quadratic_1}. Relying on the existence and uniqueness result of \cite{kpz_quadratic_1}, we define the corresponding $g$-expectations and study some of their properties. We obtain in particular a non-linear Doob-Meyer decomposition for $g$-submartingales and a downcrossing inequality which implies their regularity in time. As a consequence of these results, we also obtain a converse comparison theorem for our class of BSDEs. Finally, we provide a dual representation for the corresponding dynamic risk measures, and study the properties of their inf-convolution, giving several explicit examples.

\vspace{0.8em}
\noindent{\bf Key words:} BSDEs, quadratic growth, jumps, non-linear Doob-Meyer decomposition, dynamic risk measures, inf-convolution.
\vspace{5mm}

\noindent{\bf AMS 2000 subject classifications:} 60H10, 60H30
\end{abstract}
\newpage

\section{Introduction}

Motivated by duality methods and maximum principles for optimal stochastic control, Bismut studied in \cite{bis} a linear backward stochastic differential equation (BSDE). In their seminal paper \cite{pardpeng}, Pardoux and Peng generalized such equations to the non-linear Lipschitz case and proved existence and uniqueness results in a Brownian framework. Since then, a lot of attention has been given to BSDEs and their applications, not only in stochastic control, but also in theoretical economics, stochastic differential games and financial mathematics. 

\vspace{0.5em}
Let us now precise the structure of these equations in a discontinuous setting. Given a filtered probability space $(\Omega,\mathcal F,\left\{\mathcal F_t\right\}_{0\leq t\leq T},\mathbb P)$ generated by an $\mathbb R^d$-valued Brownian motion $B$ and a random measure $\mu$ with compensator $\nu$, solving a BSDEJ with generator $g$ and terminal condition $\xi$ consists in finding a triple of progressively measurable processes $(Y,Z,U)$ such that for all $t \in [0,T]$, $\mathbb P-a.s.$
\begin{align}
 Y_t=\xi +\int_t^T g_s(Y_s,Z_s,U_s)ds-\int_t^T Z_s dB_s -\int_t^T \int_{\R^d\backslash \{0\}} U_s(x)(\mu-\nu)(ds,dx). \label{def_bsdej}
\end{align}
We refer the reader to Section \ref{notations_qbsdej} for more precise definitions and notations.

\vspace{0.5em}
In this paper, $g$ will be supposed to satisfy a Lipschitz-quadratic growth property. More precisely, $g$ will be Lipschitz in $y$, and will satisfy a quadratic growth condition in $(z,u)$ (see Assumption \ref{assump:hquad}(iii) below). The interest for such a class of quadratic BSDEs has increased a lot in the past few years, mainly due to the fact that they naturally appear in many stochastic control problems, for instance involving utility maximization (see among many others \cite{ekr} and \cite{him}). When the filtration is generated only by a Brownian motion, the existence and uniqueness of quadratic BSDEs with a bounded terminal condition has been first treated by Kobylanski \cite{kob}. Then Tevzadze \cite{tev} introduced a new approach, consisting of a direct proof in the Lipschitz-quadratic setting. He uses a fixed-point argument to obtain existence of a solution for small terminal condition, and then pastes solutions together in the general bounded case. We refer the reader to our paper \cite{kpz_quadratic_1} for more references on the class of quadratic BSDEs.

\vspace{0.5em}
In our accompanying paper \cite{kpz_quadratic_1}, we extended the fixed-point methodology of Tevzadze to the case of a discontinuous filtration. We proved an existence and uniqueness result for bounded solutions of quadratic BSDEs. We used a comparison theorem to deduce our uniqueness result. Nonetheless, in this framework with jumps, we need additional assumptions on the generator $g$ for a comparison theorem to hold. We used either the Assumption \ref{assump.roy}, first introduced by Royer \cite{roy}, or a convexity assumption on $g$, which was already considered by Briand and Hu \cite{bh2} in the continuous case.

\vspace{0.5em}
This wellposedness result for bounded quadratic BSDEs with jumps opens the way to many possible applications. We can consider the solution of a BSDE as an operator acting on the terminal condition, this is the point of view of the $g$-expectations. It has been introduced by Peng \cite{peng1} as an example of non-linear expectation. The $g$-expectations have been extended to the case of quadratic coefficients by Ma and Yao \cite{ma}, or to discontinuous filtrations by Royer \cite{roy} and Lin \cite{lin}. It is natural in this context to use these non-linear expectations to define non-linear sub- and supermartingales (see Definition \ref{def_gsubmartingale}). In this paper, we go further in the study of quadratic BSDEs with jumps by proving a non-linear Doob Meyer decomposition for $g$-submartingales. As a consequence, we also obtain a converse comparison theorem. These results hold true under the same assumptions as the ones needed for the comparison theorem.

\vspace{0.5em}
When the generator is convex, we obtain a convex operator, which is then naturally used to construct examples of dynamic convex risk measures. Barrieu and El Karoui \cite{elkarbar2} used quadratic BSDEs to define time consistent convex risk measures and study their properties. We extend here some of these results to the case with jumps. We prove an explicit dual representation of the solution $Y$, when $g$ is independent of $y$ and convex in $(z,u)$. This allows to study some particular risk measures on a discontinuous filtration, like the \textit{entropic risk measure}, corresponding to the solution of a quadratic BSDE. Finally, we prove an explicit representation for the inf-convolution of quadratic BSDEs, thus giving the form of the optimal risk transfer between two agents using quadratic convex $g$-expectations as risk measures. The inf-convolution is again a convex operator, solving a particular BSDE. We give a sufficient condition for this BSDE to have a coefficient satisfying a quadratic growth property.

\vspace{0.5em}
The rest of this paper is organized as follows. In Section \ref{section.1}, we recall the notations, assumptions and main results of \cite{kpz_quadratic_1}, then in Section \ref{section.g}, we study general properties of quadratic $g$-martingales with jumps, such as regularity in time and the Doob-Meyer decomposition. Finally, Section \ref{section.app} is devoted to the analysis of a dual representation of the corresponding dynamic convex risk measures and to the calculation of their inf-convolution.

\section{Preliminaries} \label{section.1}

%
%
We consider in all the paper a filtered probability space $\left(\Omega,\mathcal F, \left\lbrace \mathcal F_t\right\rbrace_{0\leq t\leq T},\mathbb P\right)$, whose filtration satisfies the usual hypotheses of completeness and right-continuity. We suppose that this filtration is generated by a $d$-dimensional Brownian motion $B$ and an independent integer valued random measure $\mu(\omega,dt,dx)$ defined on $\mathbb R^+\times E$, with compensator $\lambda(\omega,dt,dx)$. $\widetilde \Omega:= \Omega \times \mathbb R^+ \times E$ is equipped with the $\sigma$-field  $\widetilde{\mathcal P}:= \mathcal P \times \mathcal E$, where $\mathcal P$ denotes  the predictable $\sigma$-field on $\Omega \times \mathbb R^+$ and $\Ec$ is the Borel $\sigma$-field on $E$.

\vspace{0.5em}
To guarantee the existence of the compensator $\lambda(\omega,dt,dx)$, we assume that for each $A$ in $\Bc(E)$ and each $\omega$ in $\Omega$, the process $X_t:= \mu(\omega,A,[0,t]) \in \Ac^+_{loc}$, which means that there exists an increasing sequence of stopping times $(T_n)$ such that $T_n \to + \infty$ a.s. and the stopped processes $X_t^{T_n}$ are increasing, c\`adl\`ag, adapted and satisfy $\E[X_{\infty}]<+\infty$.

\vspace{0.5em}
We assume in all the paper that $\lambda$ is absolutely continuous with respect to the Lebesgue measure $dt$, i.e. $\lambda(\omega,dt,dx)=\nu_t(\omega,dx)dt$. Finally, we denote $\widetilde\mu$ the compensated jump measure
$$\widetilde\mu(\omega,dx,dt) = \mu(\omega,dx,dt) - \nu_t(\omega,dx)\, dt.$$

We introduce for $1<p\leq +\infty$ the spaces
$$L^p(\nu):=\left\{u,\ \text{$\mathcal E$-measurable, such that $u\in L^p(\nu_t)$ for all $0\leq t\leq T$}\right\}.$$

Since the compensator $\nu$ depends on $\omega$, the martingale representation property do not necessarily hold. That is why we make the following assumption.

\begin{Assumption}\label{martingale_representation}
Any local martingale $M$ with respect to the filtration $(\mathcal F_t)_{0\leq t\leq T}$ has the predictable representation property, that is to say that there exist a unique predictable process $H$ and a unique predictable function $U$ such that $(H,U)\in\mathcal Z\times\mathcal U$ and
$$M_t=M_0+\int_0^tH_sdB_s+\int_0^t\int_EU_s(x)\widetilde\mu(dx,ds), \; \mathbb P-a.s.$$
\end{Assumption}

\begin{Remark}
This martingale representation property holds for instance when the compensator $\nu$ does not depend on $\omega$, i.e when $\nu$ is the compensator of the counting measure of an additive process in the sense of Sato \cite{sato}. It also holds when $\nu$ has the particular form described in \cite{kpz1}, in which case $\nu$ depends on $\omega$.
\end{Remark} 

\subsection{Notations} \label{notations_qbsdej}
We introduce the following norms and spaces for any $p\geq 1$.

\vspace{0.5em}
$\mathcal S^p$ is the space of $\mathbb R$-valued c\`adl\`ag and $(\mathcal F_t)$-progressively measurable processes $Y$ such that
$$\No{Y}^p_{\mathcal S^p}:=\mathbb E\left[\underset{0\leq t\leq T}{\sup}Y_t^p\right]<+\infty.$$

$\mathcal S^\infty$ is the space of $\mathbb R$-valued c\`adl\`ag and $(\mathcal F_t)$-progressively measurable processes $Y$ such that
$$\No{Y}_{\mathcal S^\infty}:=\underset{0\leq t\leq T}{\sup}\No{Y_t}_\infty<+\infty.$$

$\mathbb H^p$ is the space of $\mathbb R^d$-valued and $(\mathcal F_t)$-predictable processes $Z$ such that
$$\No{Z}^p_{\mathbb H^p}:=\mathbb E\left[\left(\int_0^T\abs{Z_t}^2dt\right)^{\frac p2}\right]<+\infty.$$


$\mathbb J^p$ is the space of predictable and $\mathcal E$-measurable applications $U:\Omega\times[0,T]\times E$ such that
$$\No{U}^p_{\mathbb J^p}:=\mathbb E\left[\left(\int_0^T\int_E\abs{U_s(x)}^2\nu_s(dx)ds\right)^{\frac p2}\right]<+\infty.$$

\vspace{0.5em}
Following Tang and Li \cite{li} and Barles et al. \cite{barles}, the definition of a BSDE with jumps is then 

\begin{Definition}\label{def_bsdej2}
Let $\xi$ be a $\mathcal F_T$-measurable random variable. A solution to the BSDEJ with terminal condition $\xi$ and generator $g$ is a triple $(Y,Z,U)\in\mathcal S^2\times\mathbb H^2\times\mathbb J^2 $ such that
\begin{equation}
Y_t=\xi+\int_t^Tg_s(Y_s,Z_s,U_s)ds-\int_t^TZ_sdB_s-\int_t^T\int_{E} U_s(x)\widetilde\mu(dx,ds),\ 0\leq t\leq T,\ \mathbb P-a.s.
\label{eq:bsdej}
\end{equation}
\end{Definition}
where $g:\Omega\times[0,T]\times\mathbb R\times\mathbb R^d\times \Ac(E) \rightarrow \mathbb R$ is a given application and 
$$\Ac(E):=\left\{u: \, E \rightarrow \R,\, \Bc(E)-\text{measurable} \right\}.$$



For later use, we also introduce the following BMO-type spaces. $\rm{BMO}$ is the space of square integrable c\`adl\`ag $\mathbb R^d$-valued martingales $M$ such that
$$\No{M}_{\rm{BMO}}:=\underset{\tau\in\mathcal T_0^T}{\esup^\P}\No{\mathbb E_\tau\left[\left(M_T-M_{\tau^-}\right)^2\right]}_{\infty}<+\infty,$$
where for any $t\in[0,T]$, $\mathcal T_t^T$ is the set of $(\mathcal F_s)_{0\leq s\leq T}$-stopping times taking their values in $[t,T]$.

\vspace{0.3em}
$\mathbb J^2_{\rm{BMO}}$ is the space of predictable and $\mathcal E$-measurable applications $U:\Omega\times[0,T]\times E$ such that
$$\No{U}^2_{\mathbb J^2_{\rm{BMO}}}:=\No{\int_0^.\int_EU_s(x)\widetilde\mu(dx,ds)}_{ \rm{BMO}}<+\infty.$$

\vspace{0.3em}
$\mathbb H^2_{\rm{BMO}}$ is the space of $\mathbb R^d$-valued and $\mathcal F_t$-progressively measurable processes $Z$ such that
$$\No{Z}^2_{\mathbb H^2_{\rm{BMO}}}:=\No{\int_0^.Z_sdB_s}_{\rm{BMO}}<+\infty.$$

\subsection{The non-linear generator}

We now give our quadratic growth assumption on the generator $g$. 

\begin{Assumption}\label{assump:hquad}[Quadratic growth]

\vspace{0.5em}
\rm{(i)} For fixed $(y,z,u)$, $g$ is $\mathbb{F}$-progressively measurable.

\vspace{0.5em}
\rm{(ii)} For any $p\geq 1$
\begin{equation}\label{inte}
\underset{\tau\in\mathcal T_0^T}{\esup^\P}\ \mathbb E_\tau\left[\left(\int_\tau^T\abs{g_t(0,0,0)}dt\right)^p\right]<+\infty, \ \mathbb P-a.s.
\end{equation} 

\rm{(iii)} $g$ has the following growth property. There exist $(\beta,\gamma)\in \mathbb R_+\times \mathbb R^*_+$ and a positive predictable process $\alpha$ satisfying the same integrability condition \reff{inte} as $g_t(0,0,0)$, such that for all $(\omega,t,y,z,u)$
\begin{align*}
 -\alpha_t-\beta\abs{y}-\frac\gamma2\abs{z}^2-\frac{j_t(-\gamma u)}{\gamma} \leq g_t(\omega,y,z,u)-g_t(0,0,0)\leq \alpha_t+\beta\abs{y}+\frac\gamma2\abs{z}^2+\frac{j_t(\gamma u)}{\gamma},
\end{align*}
where $j_t(u):=\int_E\left(e^{u(x)}-1-u(x)\right)\nu_t(dx).$ 
\end{Assumption}

%
Notice that $j$ is well defined on $L^2(\nu)\cap L^\infty(\nu)$. The next assumption is needed for our existence result to hold. It concerns the regularity in the $y$ variable and the differentiability in $z$ and $u$.

\begin{Assumption}\label{assump:hh}
\begin{itemize}
\item[\rm{(i)}] $g$ is uniformly Lipschitz in $y$.
$$\abs{g_t(\omega,y,z,u)-g_t(\omega,y',z,u)}\leq C\abs{y-y'}\text{ for all }(\omega,t,y,y',z,u).$$
\item[\rm{(ii)}] $g$ is $C^2$ in $z$ and there is $\theta>0$ and a process $(r_t)_{0\leq t\leq T}\in\mathbb H^2_{\rm{BMO}}$,  s.t. for all $(t,\omega,y,z,u)$,
$$\lvert D_z g_t(\omega,y,z,u)\rvert\leq r_t + \theta\abs{z}, \ \lvert D^2_{zz} g_t(\omega,y,z,u)\rvert\leq \theta.$$
\item[\rm{(iii)}] $g$ is twice Fr\'echet differentiable in the Banach space $L^2(\nu)$ and there are constants $\theta$, $\delta>0$, $C_1\geq-1+\delta$, $C_2\geq 0$ and a predictable function $m\in\mathbb J^2_{\rm{BMO}}$ s.t. for all $(t,\omega,y,z,u,x)$,
$$\abs{D_u g_t(\omega,y,z,u)}\leq m_t+\theta\abs{u},\ C_1(1\wedge\abs{x})\leq D_ug_t(\omega,y,z,u)(x)\leq C_2(1\wedge\abs{x})$$
$$ \No{D^2_u g_t(\omega,y,z,u)}_{L^2(\nu_t)}\leq \theta.$$
\end{itemize}
\end{Assumption}

\begin{Remark}\label{rem_assump_lip}
The assumption $(i)$ above is classic in the BSDE theory. The assumptions $(ii)$ and $(iii)$ are generalizations to the jump case of the assumptions considered by Tevzadze \cite{tev}. They are useful in our proof of existence in \cite{kpz_quadratic_1}. Moreover, we recall that Assumption \ref{assump:hh} implies the following,

\vspace{0.3em}
$\bullet$ There exists $\mu>0$ such that for all $(t,y,z,z',u)$
$$\abs{ g_t(\omega,y,z,u)- g_t(\omega,y,z',u)-\phi_t.(z-z')}\leq \mu \abs{z-z'}\left(\abs{z}+\abs{z'}\right),$$
where $\phi_t:=D_zg_t(y,0,u)\in\mathbb H^2_{\rm{BMO}}$.

\vspace{0.3em}
$\bullet$ Analogously, there exists $\mu>0$ such that for all $(\omega,t,y,z,u,u')$
$$\abs{ g_t(\omega,y,z,u)- g_t(\omega,y,z,u')-\langle \psi_t,u-u'\rangle_t}\leq \mu \No{u-u'}_{L^2(\nu_t)}\left(\No{u}_{L^2(\nu_t)}+\No{u'}_{L^2(\nu_t)}\right),$$
where $\psi_t:=D_ug_t(y,z,0)\in\mathbb J^2_{\rm{BMO}}$.
\end{Remark}

Finally, in order to have a comparison theorem, we need to impose either one of the following hypothesis. The first one has been first introduced by Royer \cite{roy}, it implies that the generator $g$ is Lipschitz in $u$. The second one is a convexity assumption, it has the advantage of keeping the generator quadratic in $u$.

\begin{Assumption}\label{assump.roy}
For every $(y,z,u,u')$ there exists a predictable and $\mathcal E$-measurable process $(\gamma_t)$ such that
$$g_t(y,z,u)-g_t(y,z,u')\leq\int_E\gamma_t(x)(u-u')(x)\nu_t(dx),$$
where there exist constants $C_2>0$ and $C_1\geq-1+\delta$ for some $\delta>0$ such that
$$C_1(1\wedge\abs{x})\leq\gamma_t(x)\leq C_2(1\wedge\abs{x}).$$
\end{Assumption}

\begin{Assumption}\label{assump.bh}
$g$ is jointly convex in $(z,u)$.
\end{Assumption}

We proved in \cite{kpz_quadratic_1} the following comparison theorem, used to derive our uniqueness result.

\begin{Proposition}\label{prop.comp}
Let $\xi^1$ and $\xi^2$ be two $\mathcal F_T$-measurable random variables. Let $g^1$ be a function satisfying either of the following
\begin{enumerate}
	\item[\rm{(i)}] Assumptions \ref{assump:hquad}, \ref{assump:hh}(i),(ii) and \ref{assump.roy}.
	\item[\rm{(ii)}] Assumptions \ref{assump:hquad}, \ref{assump:hh}(i)  and \ref{assump.bh}, and that $\abs{g^1(0,0,0)}+\alpha\leq M$ where $\alpha$ is the process appearing in Assumption \ref{assump:hquad}(iii) and $M$ is a positive constant.
\end{enumerate}

Let $g^2$ be another function and for $i=1,2$, let $(Y^i,Z^i,U^i)$ be the solution of the BSDEJ with terminal condition $\xi^i$ and generator $g^i$ (we assume that existence holds in our spaces), that is to say for every $t\in[0,T]$
$$Y^i_t=\xi^i+\int_t^Tg^i_s(Y^i_s,Z^i_s,U^i_s)ds-\int_t^TZ^i_sdBs-\int_t^T\int_EU^i_s(x)\widetilde\mu(dx,ds),\ \mathbb P-a.s.$$
Assume further that $\xi^1\leq\xi^2,\ \mathbb P-a.s.$ and $g_t^1(Y_t^2,Z_t^2,U_t^2)\leq g_t^2(Y_t^2,Z_t^2,U_t^2),\ \mathbb P-a.s.$ Then $Y_t^1\leq Y_t^2$, $\mathbb P-a.s.$ Moreover in case (i), if in addition we have $Y^1_0=Y^2_0$, then for all $t$, $Y^1_t=Y^2_t$, $Z_t^1=Z_t^2$ and $U_t^1=U_t^2$, $\mathbb P-a.s.$
\end{Proposition}
%
%
The following is our main existence and uniqueness result, stated in \cite{kpz_quadratic_1}.

\begin{Theorem}\label{th.unique}
Assume that $\xi\in\mathbb L^\infty$, and that the generator $g$ satisfies either 
\begin{enumerate}
	\item[\rm{(i)}] Assumptions \ref{assump:hquad}, \ref{assump:hh}(i),(ii) and \ref{assump.roy}.
	\item[\rm{(ii)}] Assumptions \ref{assump:hquad}, \ref{assump:hh} and \ref{assump.bh}, and that $g(0,0,0)$ and the process $\alpha$ appearing in Assumption \ref{assump:hquad}(iii) are bounded by some constant $M>0$.
\end{enumerate}
Then there exists a unique solution to the BSDEJ \reff{eq:bsdej}.
\end{Theorem}

\section{Quadratic non-linear expectations with jumps}\label{section.g}

The theory of $g$-expectations was introduced by Peng in \cite{peng1} as an example of non-linear expectations. Since then, numerous authors have generalized his results, extending them notably to the case of quadratic coefficients (see Ma and Yao \cite{ma}). An extension to discontinuous filtrations was obtained by Royer \cite{roy} and Lin \cite{lin}. In particular, Royer \cite{roy} gave domination conditions under which a non-linear expectation is a $g$-expectation. We refer the interested reader to these papers for more details, and we recall for simplicity some of their general properties below. Let us start with a general definition.
\begin{Definition}\label{def_gsubmartingale}
Let $\xi\in\mathbb L^\infty$ and let $g$ be such that the BSDEJ with generator $g$ and terminal condition $\xi$ has a unique solution and such that comparison in the sense of Proposition \ref{prop.comp} holds (for instance $g$ could satisfy any of the conditions in Theorem \ref{th.unique}). Then for every $t\in[0,T]$, we define the conditional $g$-expectation of $\xi$ as follows
$$\mathcal E^g_t[\xi]:=Y_t,$$
where $(Y,Z,U)$ solves the following BSDEJ
$$Y_t=\xi+\int_t^Tg_s(Y_s,Z_s,U_s)ds-\int_t^TZ_sdB_s-\int_t^T\int_EU_s(x)\widetilde\mu(dx,ds).$$
\end{Definition}

\begin{Remark}
Notice that $\Ec^g: \, \L^{\infty}(\Omega, \Fc_T, \P) \rightarrow \L^{\infty}(\Omega, \Fc_t, \P)$ does not define a true operator. Indeed, to each bounded $\Fc_T$-measurable random variable $\xi$, we associate the value $Y_t$, which is defined $\P$-a.s., i.e. outside a $\P$-negligible set $N$, but this set $N$ depends on $\xi$. We cannot {\it {\it a priori}} find a common negligible set for all variables in $\L^{\infty}$, and then define an operator $\Ec^g$ on a fixed domain, except if we only consider a countable set of variables $\xi$ on which acts $\Ec^g$.
\end{Remark}

There is a notion of g-martingales and g-sub(super)martingales.
\begin{Definition}
$X\in\mathcal S^\infty$ is called a $g$-submartingale (resp. $g$-supermartingale) if 
$$\mathcal E^g_s[X_t]\geq\text{ (resp. $\leq$)}X_s,\ \mathbb P-a.s.,\text{ for any $0\leq s\leq t\leq T$}.$$

$X$ is called a $g$-martingale if it is both a $g$-sub and supermartingale.
\end{Definition}

The following results are easy generalizations of the classical arguments which can be found in \cite{peng1} or \cite{elkarbar2}, and are consequences of the comparison theorem. We therefore omit the proofs.
\begin{Lemma}
 $\left\{ \Ec^g_t \right\}_{t\geq 0}$ is monotonic increasing and time consistent, i.e.
  \begin{itemize}
   \item $\xi_1 \geq \xi_2$, $\P$-a.s. implies that $\Ec^g_t(\xi_1) \geq \Ec^g_t(\xi_2)$, $\P$-a.s., $\forall t \geq 0$.
    \item For any bounded stopping times $R \leq S \leq \tau$ and $\Fc_{\tau}$-measurable random variable $\xi_{\tau}$,
\begin{align}
 \Ec^g_R(\Ec^g_S(\xi_{\tau})) = \Ec^g_R(\xi_{\tau}) \; \P\text{-a.s.}
\end{align}
  \end{itemize}
\end{Lemma}

\begin{Definition}
 We will say that $\Ec^g$ is
 
\vspace{0.3em}
(i) Constant additive, if for any stopping times $R \leq S$, any $\Fc_R$-measurable random variable $\eta_R$ and any $\Fc_S$-measurable random variable $\xi_S$, $$\Ec^g_R(\xi_S + \eta_R) = \Ec^g_R(\xi_S) + \eta_R, \ \P \text{-a.s.}$$

\vspace{0.3em}
(ii) Positively homogeneous, if for any stopping times $R \leq S$, and any positive $\Fc_R$-measurable random variable $\lambda$, $$\Ec^g_R(\lambda \xi_S) = \lambda \Ec^g_R(\xi_S).$$

\vspace{0.3em}
(iii) Convex, if for any stopping times $R \leq S$, any random variables $(\xi^1_S, \xi^2_S)$ and any $\lambda \in [0,1]$,
	  $$ \Ec^g_R(\lambda \xi^1_S + (1-\lambda)\xi^2_S) \leq \lambda \Ec^g_R(\xi^1_S) + (1-\lambda) \Ec^g_R(\xi^2_S).$$
\end{Definition}

\vspace{0.5em}
The next Lemma shows that the operator $\mathcal E^g$ inherits the above properties from $g$.
\begin{Lemma}\label{g.prop.Eg}
(i) If $g$ does not depend on $y$, then $\Ec^g$ is constant additive.

\vspace{0.3em}
(ii) If $g$ is positively homogeneous in $(y,z,u)$, then $\Ec^g$ is positively homogeneous.

\vspace{0.3em}
(iii) If $g$ is moreover right continuous on $[0,T)$ and continuous at $T$, then the reverse implications of $(i)$ and $(ii)$ are also true.

\vspace{0.3em}
(iv) $\Ec^g$ is convex if $g$ is convex in $(y,z,u)$.

\vspace{0.3em}
(v) If $g^1 \leq g^2$, $\P$-a.s., then $\Ec^{g^1} \leq \Ec^{g^2}$. If $g^1$ and $g^2$ are moreover right continuous on $[0,T)$ and continuous at $T$, then the reverse is also true.\end{Lemma}

\proof We adapt the ideas of the proofs in \cite{elkarbar2} to our context with jumps.

\vspace{0.5em}
$\rm{(i)}$ The proof of the first property is exactly the same as the proof of Theorem 6.7.b2 in \cite{elkarbar2}, so we omit it.

\vspace{0.5em}
$\rm{(ii)}$ Let $g^{\lambda}(t,y,z,u) := \frac{1}{\lambda}g_t( \lambda y, \lambda z, \lambda u)$. Then $\left\{ \frac{1}{\lambda} \Ec^g_t(\lambda \xi_S)\right\}_{t \geq 0}$ is a solution of the BSDEJ with coefficient $g^{\lambda}$ and terminal condition $\xi_S$. If $g=g^{\lambda}$, then $ \Ec^g_t(\lambda \xi_S) = \lambda\Ec^g_t(\xi_S).$

\vspace{0.5em}
$\rm{(iii)}$ The reverse implications in (i) and (ii) are direct consequences of Corollary \ref{reverse}.

\vspace{0.5em}
$\rm{(iv)}$ Suppose that $g$ is convex in $(y,z,u)$. Let $(Y^i,Z^i,U^i)$ be the unique solution of the BSDEJ with coefficients $(g,\xi^i_S)$, $i=1,2$, and set 
  $$\Yt_t = \lambda Y^1_t + (1-\lambda)Y^2_t, \ \Zt_t = \lambda Z^1_t + (1-\lambda)Z^2_t\text{ and }\Ut_t(\cdot) = \lambda U^1_t(\cdot) + (1-\lambda)U^2_t(\cdot).$$

We have
\begin{align*}
 -d\Yt_t =& \left[\lambda g_t(Y^1_t,Z^1_t,U^1_t) + (1-\lambda)g_t(Y^2_t,Z^2_t,U^2_t)\right]dt -\left(\lambda Z^1_t + (1-\lambda)Z^2_t\right)dB_t \\
 &- \int_E (\lambda U^1_t(x) + (1-\lambda)U^2_t(x))\widetilde{\mu}(dt,dx) \\
 =& \left[g_t(\Yt_t, \Zt_t, \Ut_t) + k(t,Y^1_t,Y^2_t,Z^1_t,Z^2_t,U^1_t,U^2_t,\lambda) \right]dt - \Zt_t dB_t - \int_E \Ut_t(x)\widetilde{\mu}(dt,dx),
\end{align*}
where $$k(t,Y^1_t,Y^2_t,Z^1_t,Z^2_t,U^1_t,U^2_t,\lambda):= \lambda g_t(Y^1_t,Z^1_t,U^1_t) + (1-\lambda)g_t(Y^2_t,Z^2_t,U^2_t) - g_t(\Yt_t, \Zt_t, \Ut_t),$$ is a non negative function. Then using Proposition \ref{prop.comp} we obtain in particular
$$ \Ec^g_t(\lambda \xi^1_S + (1-\lambda)\xi^2_S) \leq \Yt_t = \lambda \Ec^g_t(\xi^1_S) + (1-\lambda) \Ec^g_t(\xi^2_S).$$

\vspace{0.5em}
$\rm{(v)}$ This last property is a direct consequence of the comparison Theorem \ref{prop.comp}. The reverse implication is again a consequence of Corollary \ref{reverse}.
\ep

\begin{Example}\label{backward.riskmeasures}
 These easy properties allow us to construct examples of time consistent dynamic convex risk measures, by appropriate choices of generator $g$. 

\vspace{0.3em}
$\bullet$ Defining $g_t(z,u) := \frac\gamma2\abs{z_t}^2+\frac1\gamma j_t(\gamma u_t)$, we obtain the so called entropic risk measure on our particular filtration.

\vspace{0.3em}
$\bullet$ As proved in \cite{roy}, if we define $$g_t(z,u):=\eta \abs{z} + \eta \int_E (1\wedge \abs{x})u^+(x)\nu_t(dx) - C_1\int_E (1\wedge \abs{x})u^-(x)\nu_t(dx),$$ where $\eta>0$ and $-1<C_1 \leq 0$, then $\Ec^g$ is a convex risk measure with the following representation $\Ec_0^g(\xi)=\underset{\mathbb Q \in \mathbf{Q}}{\sup} \E^\mathbb Q\left[\xi\right]$,
 with 
 \begin{align*}
 \mathbf{Q}:= \Big\{& \mathbb Q, \, \frac{d\mathbb Q}{d\P}|\Fc_t=\Ec\left(\int_0^t \mu_s dB_s +\int_0^t \int_E v_s(x)\widetilde{\mu}(ds,dx)\right)\\
 										&\textrm{with $\mu$ and $v$ predictable, } \abs{\mu_s} \leq \eta, \, v^+_s(x)\leq \eta(1\wedge x), \, v^-_s(x)\leq C_1(1\wedge x) \Big\}.
 \end{align*}

\vspace{0.3em}
$\bullet$ If we define a linear generator $g$ by
 $$g_t(z,u):=\alpha z+\beta\int_E (1 \wedge \abs{x})u(x)\nu_t(dx), \, \alpha \in \R, \, \beta\geq -1+\delta \text{ for some } \delta>0,$$
 then we obtain a linear risk measure, since $\Ec^g$ will only consist of a linear expectation with respect to the probability measure $\Q$, whose Radon-Nikodym derivative is equal to $$\frac{d\Q}{d\P} = \Ec\left(\alpha B_t + \int_0^t \int_E \beta (1 \wedge \abs{x}) \widetilde{\mu}(ds,dx)\right).$$
\end{Example}

In the rest of this section, we will provide important properties of quadratic $g$-expectations and the associated $g$-martingales in discontinuous filtrations, which generalize the known results in simpler cases.

\subsection{Non-linear Doob Meyer decomposition}
We start by proving that the non-linear Doob Meyer decomposition first proved by Peng in \cite{peng2} still holds in our context. We have two different sets of assumptions under which this result holds, and they are both related to the assumptions under which our comparison theorem \ref{prop.comp} holds. From a technical point of view, our proof consists in approximating our generator by a sequence of Lipschitz generators. However, the novelty here is that because of the dependence of the generator in $u$, we cannot use the classical exponential transformation and then use some truncation arguments, as in \cite{kob} and \cite{ma}. Indeed, since $u$ lives in an infinite dimensional space, those truncation type arguments no longer work {\it a priori}. Instead, inspired by \cite{elkarbar}, we will only use regularizations by inf-convolution, which are known to work in any Banach space. 

\begin{Theorem}\label{dm.decomp}
Let $Y$ be a c\`adl\`ag $g$-submartingale (resp. $g$-supermartingale) in $\mathcal S^\infty$ (we assume that existence and uniqueness for the BSDEJ with generator $g$ hold for any bounded terminal condition). Assume further either one of these conditions 
\begin{itemize}
	\item[\rm{(i)}] Assumptions \ref{assump:hquad} and \ref{assump.roy} hold, with the addition that the process $\gamma$ does not depend on $(y,z)$ and that $\abs{g(0,0,0)}+\alpha\leq M$, where $\alpha$ is the process appearing in Assumption \ref{assump:hquad}(iii) and $M>0$ is constant. 
	\item[\rm{(ii)}] Assumptions \ref{assump:hquad}, \ref{assump:hh}(i) hold, $g$ is concave (resp. convex) in $(z,u)$, $\abs{g(0,0,0)}+\alpha\leq M$, where $\alpha$ is the process appearing in Assumption \ref{assump:hquad}(iii) and $M>0$ is constant.
	\end{itemize}
Then there exists a predictable non-decreasing (resp. non-increasing) process $A$ null at $0$ and processes $(Z,U)\in\mathbb H^2\times\mathbb J^2$ such that
$$Y_t=Y_T+\int_t^Tg_s(Y_s,Z_s,U_s)ds-\int_t^TZ_sdB_s-\int_t^T\int_EU_s(x)\widetilde\mu(dx,ds)-A_T+A_t,\ t\in[0,T].$$

\end{Theorem}

\begin{Remark}
We emphasize that the two assumptions in the above theorem are not of the same type. Indeed, Assumption \ref{assump.roy} implies that the generator $g$ is uniformly Lipschitz in $u$, which is a bit disappointing if we want to work in a quadratic context. This is why we also considered the convexity hypothesis on $g$, which allows us to retrieve a generator which is quadratic in both $(z,u)$. We do not know whether those two assumptions are necessary or not to obtain the result, but we remind the reader that our theorem encompasses the case of the so-called entropic generator, which has quadratic growth and is convex in $(z,u)$. To the best of our knowledge, this particular case which was already proved in \cite{ngou}, was the only result available in the literature up until now.
\end{Remark}
\proof First, if $Y$ is $g$-supermartingale, then $-Y$ is a $g^-$-submartingale where $g^-_t(y,z,u):=-g_t(-y,-z,-u).$ Since $g^-$ satisfies exactly the same Assumptions as $g$, and given that $g^-$ is convex when $g$ is concave, it is clear that we can without loss of generality restrict ourselves to the case of $g$-submartingales. We start with the first result.

\vspace{0.5em}
{\bf Step $1$:} Assumptions \ref{assump:hquad} and \ref{assump.roy} hold. We will approximate the generator $g$ by a sequence of functions $(g^n)$ which are uniformly Lipschitz in $(y,z)$ (recall that under the assumed assumptions, $g$ is already Lipschitz in $u$). We emphasize that unlike most of the literature on quadratic BSDEs, with the notable exception of \cite{elkarbar}, we will not use any exponential change in our proof. Building upon the results of Lepeltier and San Martin \cite{lsm}, we would like to use a sup-convolution to regularize our generator. However, due to the quadratic growth assumption in $z$, such a sup-convolution is not always well defined. Therefore, we will first use a truncation argument to bound our generator from above by a function with linear growth. Let us thus define for all $n\geq 0$
$$\widetilde{g}^n_t(y,z,u):=g_t(y,z,u)\wedge \left(M+n\abs{z}-\frac{\gamma}{2}\abs{z}^2\right),$$
where the constants $(\alpha,\gamma)$ are the ones appearing in Assumption \ref{assump:hh}(ii). It is clear that we have the following estimates
$$-M-\beta\abs{y}-\frac{\gamma}{2}\abs{z}^2-\frac1\gamma j_t\left(-\gamma u\right)\leq \widetilde g^n_t(y,z,u)\leq M+\beta\abs{y}+n\abs{z}+\frac1\gamma j_t\left(\gamma u\right),$$
and that $\widetilde g^n$ decreases pointwise to $g$. We now define for all $p\geq n\vee\beta$
$$\widetilde g^{n,p}_t(y,z,u):=\underset{(w,v)\in\mathbb Q^{d+1}}{\sup}\left\{\widetilde g_t^n(w,v,u)-p\abs{y-w}-p\abs{z-v}\right\}.$$

This function is indeed well-defined, since we have for $p\geq n$
\begin{align*}
\widetilde g^{n,p}_t(y,z,u)&\leq M+\frac1\gamma j_t(\gamma u)+\underset{(w,v)\in\mathbb Q^{d+1}}{\sup}\left\{\beta\abs{w}+n\abs{v}-p\abs{y-w}-p\abs{z-v}\right\}\\
&=M+\beta\abs{y}+n\abs{z}+\frac1\gamma j_t(\gamma u).
\end{align*}

Moreover, by the results of Lepeltier and San Martin \cite{lsm}, we know that $\widetilde g^{n,p}$ is uniformly Lipschitz in $(y,z)$ and that $\widetilde g^{n,p}(y,z,u)\downarrow g_t(y,z,u)$ as $n$ and $p$ go to $+\infty$. Finally, we define
$$g^n_t(y,z,u):=\widetilde g^{n,n}_t(y,z,u).$$

Then the $g^n$ are uniformly Lipschitz in $(y,z,u)$ and decrease pointwise to $g$. Now, we want somehow to use the fact that we know that the non-linear Doob-Meyer decomposition holds when the underlying generator is Lipschitz. But this was shown by Royer only when the generator also satisfies Assumption \ref{assump.roy}. Therefore, we will now verify that $g^n$ inherits Assumption \ref{assump.roy} from $g$. First of all, we show that this is true for $\widetilde g^n$. 

\vspace{0.5em}
Let $u^1,u^2\in L^\infty(\nu)\cap L^2(\nu)$ and fix some $(y,z)\in\mathbb R^{d+1}$. Then if we have 
$$g_t(y,z,u^1)\leq M+n\abs{z}-\frac{\gamma}{2}\abs{z}^2\text{ and } g_t(y,z,u^2)\leq M+n\abs{z}-\frac{\gamma}{2}\abs{z}^2,$$
then
$$\widetilde g^n_t(y,z,u^1)-\widetilde g^n_t(y,z,u^2)=g_t(y,z,u^1)-g_t(y,z,u^2),$$
and the result is clear with the same process $\gamma$ as the one for $g$. Similarly, if 
$$g_t(y,z,u^1)\geq M+n\abs{z}-\frac{\gamma}{2}\abs{z}^2\text{ and } g_t(y,z,u^2)\geq M+n\abs{z}-\frac{\gamma}{2}\abs{z}^2,$$
then
$$\widetilde g^n_t(y,z,u^1)-\widetilde g^n_t(y,z,u^2)=0,$$
and the desired result also follows by choosing the process $\gamma$ in Assumption \ref{assump.roy} to be $0$. Finally, if (the remaining case can be treated similarly)
$$g_t(y,z,u^1)\geq M+n\abs{z}-\frac{\gamma}{2}\abs{z}^2\text{ and } g_t(y,z,u^2)\leq M+n\abs{z}-\frac{\gamma}{2}\abs{z}^2,$$
then
$$\widetilde g^n_t(y,z,u^1)-\widetilde g^n_t(y,z,u^2)\leq M+n\abs{z}-\frac{\gamma}{2}\abs{z}^2-g_t(y,z,u^2)\leq g_t(y,z,u^1)-g_t(y,z,u^2),$$
and the desired result follows once more with the same process $\gamma$ as the one for $g$.

\vspace{0.5em}
Next, we show that $\widetilde g^{n,p}$ inherits Assumption \ref{assump.roy} from $\widetilde g^n$. Indeed, we have
$$\widetilde g^{n,p}_t(y,z,u^1)-\widetilde g^{n,p}_t(y,z,u^2)\leq \underset{(w,v)\in\mathbb Q^{d+1}}{\sup}\left\{\widetilde g^n_t(w,v,u^1)-\widetilde g^n_t(w,v,u^2)\right\},$$
which implies the result since the process $\gamma$ in Assumption \ref{assump.roy} does not depend on $(y,z)$.

\vspace{0.5em}
Let now $Y$ be a $g$-submartingale. We will now show that it is also a $g^n$-submartingale for all $n\geq 0$. Let now $\mathcal Y$ (resp. $\mathcal Y^n$) be the unique solution of the BSDEJ with terminal condition $Y_T$ and generator $g$ (resp. $g^n$). Since $g^n$ satisfies Assumption \ref{assump.roy} and is uniformly Lipschitz in $(y,z,u)$, we can apply the comparison theorem for Lipschitz BSDEJs (see \cite{roy}) to obtain 
$$Y_t\leq \mathcal Y_t\leq \mathcal Y^n_t,\ \mathbb P-a.s.$$

Hence $Y$ is a $g^n$-submartingale. We can therefore apply the Doob-Meyer decomposition in the Lipschitz case (see Theorem $1.1$ in Lin \cite{lin} or Theorem $4.1$ in Royer \cite{roy}) to obtain the existence of $( Z^n, U^n)\in\mathbb H^2\times\mathbb J^2$ and of a predictable non-decreasing process $ A^n$ null at $0$ such that
\begin{equation}
 Y_t=Y_T+\int_t^T g^n_t( Y_s, Z^n_s, U^n_s)ds-\int_t^T Z^n_sdB_s-\int_t^T\int_E U^n_s(x)\widetilde\mu(dx,ds)- A^n_T+ A^n_t.
\label{eq:9}
\end{equation}

Since $Y$ does not depend on $n$, the martingale part of \reff{eq:9} neither, which entails that $Z^n$ and $U^n$ are independent of $n$. We can rewrite \reff{eq:9} as
\begin{equation}
Y_t=Y_T+\int_t^T g^n_t( Y_s, Z_s, U_s)ds-\int_t^T Z_sdB_s-\int_t^T\int_E U_s(x)\widetilde\mu(dx,ds)-A^n_T+ A^n_t.
\label{eq:10}
\end{equation}

Since $ g^n$ converges pointwise to $ g$, the dominated convergence theorem implies that
$$\int_0^T\left( g^n_s( Y_s, Z_s, U_s)- g_s( Y_s, Z_s, U_s)\right)ds\rightarrow 0,\ \mathbb P-a.s.$$

Hence, it holds $\mathbb P-a.s.$ that for all $s\in[0,T]$
$$ A^n_s\rightarrow  A_s:= Y_s - Y_0+\int_0^s g_r( Y_r, Z_r, U_r)dr-\int_0^s Z_rdB_r-\int_0^s\int_E U_r(x)\widetilde\mu(dx,dr).$$

Furthermore, it is easy to see that $ A$ is still a predictable non-decreasing process null at $0$. 

\vspace{0.5em}
{\bf Step $2$:} The concave case.

\vspace{0.5em}
We have seen in the above proof that the main ingredients to obtain the desired decomposition are the comparison theorem and the non-linear Doob-Meyer decomposition in the Lipschitz case. As we have already seen in our comparison result of Proposition \ref{prop.comp}, Assumption \ref{assump.roy} plays, at least formally, the same role as the concavity/convexity assumption \ref{assump.bh}. Moreover, we show in the Appendix (see Proposition \ref{concave}) that the non-linear Doob-Meyer decomposition also holds in the Lipschitz case under Assumption \ref{assump.bh} instead of Assumption \ref{assump.roy}. We are therefore led to proceed exactly as in the previous step. Define thus
$$\widetilde g^n_t(y,z,u):=g_t(y,z,u)\wedge\left(M+n\abs{z}+n\No{u}_{L^2(\nu_t)}-\frac\gamma2\abs{z}^2-\frac{1}{\gamma}j_t(\gamma u)\right).$$

Then $\widetilde g^n$ is still concave as the minimum of two concave functions, converges pointwise to $g$ and verifies 
$$-M-\beta\abs{y}-\frac{\gamma}{2}\abs{z}^2-\frac1\gamma j_t\left(-\gamma u\right)\leq \widetilde g^n_t(y,z,u)\leq M+\beta\abs{y}+n\abs{z}+n\No{u}_{L^2(\nu_t)}.$$

Thanks to this estimate the following sup-convolution is well defined for $p\geq\beta\vee n$
$$\widetilde g^{n,p}_t(y,z,u):=\underset{(w,v,r)\in\mathbb Q^{d+1}\times L^2(\nu_t)}{\sup}\left\{\widetilde g^n_t(w,v,r)-p\abs{y-w}-p\abs{z-v}-p\No{u-r}_{L^2(\nu_t)}\right\},$$
and is still concave as the sup-convolution of concave functions.

\vspace{0.5em}
We can then finish the proof exactly as in Step $1$, using the comparison theorem of Proposition \ref{prop.comp} and the non-linear Doob-Meyer decomposition given by Proposition \ref{concave}.
\ep

\begin{Remark}
After obtaining this non-linear Doob-Meyer decomposition, it is interesting to wonder whether we can say anything about the non-decreasing process $A$ (apart from saying that it is predictable). For instance, since we are working with bounded $g$-supermartingales, we may think that $A$ can also be bounded. However, it is already known for classical supermartingales (corresponding to the case $g=0$) that this is not true. Indeed, let $X$ be a supermartingale and let $A$ be the predictable non-decreasing process appearing in its Doob-Meyer decomposition. Then, the inequality $\abs{X_t}\leq M$ for all $t$ only implies that
$$\mathbb E\left[(A_t)^p\right]\leq p!M^p,\text{ for all $p\geq 1$.}$$

Since we have $\mathbb E_t\left[X_t-X_T\right]=\mathbb E_t\left[A_T-A_t\right],$ we may then wonder if there could exists another non-decreasing process $C_t$ bounded but not necessarily adapted such that 
\begin{equation}
\mathbb E_t\left[X_t-X_T\right]=\mathbb E_t\left[C_T-C_t\right].
\label{eq:meyer}
\end{equation}

This result is then indeed true, and as shown by Meyer \cite{meyer}, if $X$ is c\`adl\`ag, positive, bounded by some constant $M$, then if we denote $\dot{X}$ the predictable projection of $X$, the non-decreasing process $C$ in \reff{eq:meyer} is given by
\begin{equation}
C_T-C_t=M\left(1-\exp\left(-\int_t^T\frac{dA_s^c}{M-\dot{X}_s}\right)\prod_{t<s\leq T}\left(1-\frac{\Delta A_s}{M-\dot{X}_s}\right)\right),
\label{eq:meyer2}
\end{equation}
where $A^c$ is the continuous part of $A$. If we now consider a $g$-supermartingale $Y$ satisfying either one of the assumptions in Theorem \ref{dm.decomp}, a simple application of It\^o's formula shows that 
$$\widetilde Y_t:=\exp\left(\gamma Y_t+\gamma Mt+\gamma\beta\int_0^t\abs{Y_s}ds\right),$$
is a bounded classical supermartingale, which therefore admits the following decomposition
\begin{equation}\label{doobm}
\widetilde Y_t=\widetilde Y_0+\int_0^t\widetilde Z_sdB_s+\int_0^t\int_E\widetilde U_s(x)\widetilde\mu(dx,ds) +\widetilde A_t,\ \mathbb P-a.s.,
\end{equation}
for some $(\widetilde Z, \widetilde U)\in\mathbb H^2\times\mathbb J^2$ and some predictable non-decreasing process $\widetilde A$. We can the apply Meyer's result to obtain
$\mathbb E_t[\widetilde Y_t-\widetilde Y_T]=\mathbb E_t\left[D_T-D_t\right],$ where $D$ is given by \reff{eq:meyer2}.

\vspace{0.5em}
Then, applying It\^o's formula to $\ln(\widetilde Y_t)$ in \reff{doobm}, we can show after some calculations that
$$Y_t=Y_T+\int_t^T g_t( Y_s, Z_s, U_s)ds-\int_t^T Z_sdB_s-\int_t^T\int_E U_s(x)\widetilde\mu(dx,ds)+A_T-A_t,$$
where $(Z,U)\in\mathbb H^2\times\mathbb J^2$ and $A$ is a predictable process with finite variation, and where $(Z,U,A)$ can be computed explicitly from $(\widetilde Z,\widetilde U,\widetilde A)$. 

\vspace{0.5em}
By uniqueness of the non-linear Doob-Meyer decomposition for $Y$, $A$ is actually non-decreasing,  and we have a result somehow similar to that of Meyer, using the relation between $\widetilde A$ and $A$. It would of course be interesting to pursue further this study.
\end{Remark}

We end this section with a converse comparison result for our class of quadratic BSDEJs, which is a consequence of the previous Doob-Meyer decomposition.

\begin{Corollary}\label{reverse}
Let $g^1$ be a function satisfying either one of the assumptions in Theorem \ref{dm.decomp} and $g^2$ be another function. We furthermore suppose that $t \mapsto g^i_t(\cdot,\cdot,\cdot)$ is right continuous in $t \in [0,T)$ and continuous at $T$, for $i=1,2$. For any $\xi\in L^\infty$, denote for $i=1,2$, $Y_t^{i,\xi}$ the solution of the BSDEJ with generator $g^i$ and terminal condition $\xi$ (existence and uniqueness are assumed to hold in our spaces). If we have 
$$Y^{1,\xi}_t \leq Y^{2,\xi}_t, \ t\in[0,T],\ \forall\xi\in L^\infty, \ \P-a.s.,$$
then we have 
$$g^1_t(y,z,u) \leq g^2_t(y,z,u), \, \forall (t,y,z,u), \, \P-a.s.$$
\end{Corollary}

\proof
For any $\xi\in L^\infty$, the assumption of the Corollary is equivalent to saying that $Y^{2,\xi}$ is a $g^1$-supermartingale. Given the assumptions on $g^1$, we can apply Theorem \ref{dm.decomp} to obtain the existence of $(\widetilde Z^{2,\xi},\widetilde U^{2,\xi},A^{2,\xi})$ such that for any $0\leq s<t\leq T$ we have, $\mathbb P-a.s.$
\begin{equation}
Y_s^{2,\xi}=Y_t^{2,\xi}+\int_s^tg_r^1\left(Y^{2,\xi}_r, \widetilde Z^{2,\xi}_r,\widetilde U_r^{2,\xi}\right)dr-\int_s^t\widetilde Z_r^{2,\xi}dB_r-\int_s^t\int_E \widetilde U_r^{2,\xi}(x)\widetilde\mu(dx,dr)+A_t^{2,\xi}-A_s^{2,\xi}.
\label{trou}
\end{equation}

Moreover, if we denote $(Y^{2,\xi},Z^{2,\xi},U^{2,\xi})$ the solution of the BSDEJ with generator $g^2$ and terminal condition $\xi$, we also have by definition
\begin{equation}
Y_s^{2,\xi}=Y_t^{2,\xi}+\int_s^tg_r^2\left(Y^{2,\xi}_r,Z^{2,\xi}_r,U_r^{2,\xi}\right)dr-\int_s^t Z_r^{2,\xi}dB_r-\int_s^t\int_E U_r^{2,\xi}(x)\widetilde\mu(dx,dr),\ \mathbb P-a.s.
\label{trou2}
\end{equation}

Identifying the martingale parts in \reff{trou} and \reff{trou2}, we obtain that $\mathbb P-a.s.$, $\widetilde Z^{2,\xi}=Z^{2,\xi}$ and $\widetilde U^{2,\xi}=U^{2,\xi}$. Furthermore, this implies by taking the expectation that
$$\frac{1}{t-s}\int_s^t\mathbb E\left[g^1_r\left(Y_r^{2,\xi},Z^{2,\xi}_r,U^{2,\xi}_r\right)\right]dr\leq \frac{1}{t-s}\int_s^t\mathbb E\left[g^2_r\left(Y^{2,\xi}_r,Z^{2,\xi}_r,U^{2,\xi}_r\right)\right]dr.$$

Now, we finish using the same argument as in Chen \cite{chen}. Let $\xi=X_T$ where for a given $(s,y_0,z_0,u_0)$, $X$ is the solution of the SDE (existence and uniqueness are classical, see for instance Jacod \cite{jacques})
$$X_t=y_0-\int_s^tg^2_r(X_r,z_0,u_0)dr+\int_s^tz_0dB_r+\int_s^t\int_Eu_0(x)\widetilde\mu(dx,dr).$$

Letting $t\longrightarrow s^+$, we obtain $g^1_s(y_0,z_0,u_0)\leq g^2_s(y_0,z_0,u_0),$ which is the desired result.
\ep

\subsection{Upcrossing inequality}

In this subsection, we prove an upcrossing inequality for quadratic $g$-submartingales, which is similar to the one obtained by Ma and Yao \cite{ma} in the case without jumps. This property is essential for the study of path regularity of $g$-submartingales. 

\begin{Theorem}\label{up}
Let $(X_t)$ be a $g$-submartingale (reps. $g$-supermartingale) and assume that either one of the following holds (as usual we assume existence and uniqueness for the solution of the BSDEJ driven by $g$ with any bounded terminal condition)
\begin{itemize}
	\item[\rm{(i)}] Assumptions \ref{assump:hquad}, \ref{assump:hh}(i),(ii) and \ref{assump.roy} hold, with the addition that $\abs{g(0,0,0)}+\alpha\leq M$, where $\alpha$ is the process appearing in Assumption \ref{assump:hquad}(iii) and $M>0$ is constant. 
	\item[\rm{(ii)}] Assumptions \ref{assump:hquad}, \ref{assump:hh}(i) hold, $g$ is concave (reap. convex), with the addition that $\abs{g(0,0,0)}+\alpha\leq M$, where $\alpha$ is the process appearing in Assumption \ref{assump:hquad}(iii) and $M>0$ is constant. 
	\end{itemize}

\vspace{0.5em}
Set $J:=\gamma M(e^{\beta T}-1)/\beta+\gamma e^{\beta T}\No{X}_{\mathcal{S}^\infty}$, 
and denote for any $\theta\in(0,1)$
$$\widetilde{X}_t:=X_t+k(J+1)t,\ \widehat{X}_t:=\exp\left(k_\theta(1+J)t+\frac{\gamma\theta}{1-\theta}X_t\right)t\in\left[0,T \right],$$
where $k$ and $k_\theta$ are a well-chosen constants depending on $\theta$, $C$, $M$, $\beta$ and $\gamma$, the constants in Assumption \ref{assump:hquad}. Let $0=t_0<t_1<...<t_n=T$ be a subdivision of $[0,T]$ and let $a<b$, 
we denote $U_a^b[\widetilde X,n]$, the number of upcrossings of the interval $[a,b]$ by $(\widetilde{X}_{t_j})_{0\leq j\leq n}$. 
Then 

\vspace{0.3em} 
$\bullet$ If (i) above holds, there exists a BMO process $(\lambda^n_t,\ t\in\left[ 0,T\right])$ such that
$$ \mathbb E\left[U_a^b[X,n]\mathcal E_T\right]\leq\frac{\No{X}_{\mathcal S^\infty}+2k(J+1)T+\abs{a}}{b-a},$$
with
$$\mathcal E_T:=\mathcal E\left(\int_0^{T}\left( \lambda^n_s +\phi_s\right) dB_s+\int_0^{T}\int_E\gamma_s(x)\widetilde\mu(dx,ds)\right),$$
and where $\phi$ and $\gamma$ are defined in Remark \ref{rem_assump_lip} and Assumption \ref{assump.roy}, and such that $$ \mathbb E\left[\int^{t_n}_0\abs{\lambda^n_s}^2ds\right]\leq C_1.$$ 

\vspace{0.3em} 
$\bullet$ If (ii) above holds, then for any $\theta\in(0,1)$
$$\mathbb E\left[U_a^b[X,n]\right]\leq\frac{\exp\left(\frac{\gamma\theta}{1-\theta}\No{X}_{\mathcal{S}^\infty}\right)+\exp\left(\frac{\gamma\theta}{1-\theta}\abs{a}\right)}{\exp\left(\frac{\gamma\theta}{1-\theta}b\right)-\exp\left(\frac{\gamma\theta}{1-\theta}a\right)}.$$
\end{Theorem}

\proof
As usual, we can restrict ourselves to the $g$-submartingale case.

\vspace{0.5em}
{\bf Step $1$:} When (i) holds. For any $j \in {1,\cdots,n}$, we consider the following BSDEJ
\begin{equation}
Y^j_t=X_{t_j}+\int_t^{t_j}g_s(Y^j_s,Z^j_s,U^j_s)ds-\int_t^{t_j}Z^j_sdB_s-\int_t^{t_j}\int_EU^j_s(x)\widetilde\mu(dx,ds),\ 0\leq t\leq {t_j},\ \mathbb P-a.s.
\label{eq:bsdej_up}
\end{equation}
From Proposition $3.1$ in \cite{kpz_quadratic_1} one has
\begin{equation}
\No{Y^j}_{\mathcal{S}^\infty}\leq \gamma M\frac{e^{\beta(t_j-t_{j-1})}-1}{\beta}+\gamma e^{\beta(t_j-t_{j-1})}\No{X_{t_j}}_{\mathcal S^\infty}\leq J.
\label{norm_estimate}
\end{equation}

We can rewrite \reff{eq:bsdej_up} as follows
\begin{align*}
Y^j_t=&\ X_{t_j}+\int_t^{t_j}\left[ g_s(Y^j_s,Z^j_s,U^j_s)-g_s(Y^j_s,0,U^j_s)\right] ds\\
&+\int_t^{t_j}\left[ g_s(Y^j_s,0,U^j_s)-g_s(0,0,U^j_s)\right] ds + \int_t^{t_j}\left[ g_s(0,0,U^j_s)-g_s(0,0,0)\right] ds\\
&+\int_t^{t_j}g_s(0,0,0) ds-\int_t^{t_j}Z^j_sdB_s-\int_t^{t_j}\int_EU^j_s(x)\widetilde\mu(dx,ds),\ 0\leq t\leq {t_j},\ \mathbb P-a.s.
\end{align*}
Then by Remark \ref{rem_assump_lip}, there exist a bounded process $\eta^n$ and $(\phi, \lambda^n)\in\mathbb H_{\rm{BMO}}^2$ with 
$$\abs{\lambda^n_t}\leq \mu\abs{Z^j_t}, \text{ }\mathbb P-a.s.,\forall t\in\left[t_{j-1},t_j \right],$$
such that
\begin{align*}
Y^j_t=&\ X_{t_j}+\int_t^{t_j}\left[\left( \lambda^n_s +\phi_s\right)Z^j_s+ \eta^n_s Y^j_s\right] ds+\int_t^{t_j}\left[ g_s(0,0,U^j_s)-g_s(0,0,0)\right] ds\\
&+\int_t^{t_j}g_s(0,0,0) ds-\int_t^{t_j}Z^j_sdB_s-\int_t^{t_j}\int_EU^j_s(x)\widetilde\mu(dx,ds)\\
\leq & \ X_{t_j}+k(J+1)(t_j-t)-\int_t^{t_j}Z^j_sdB^n_s-\int_t^{t_j}\int_EU^j_s(x) \widetilde\mu_1(ds,dx),
\end{align*}
for some positive constant $k$ and where
$$B^n_t:=B_t-\int_0^t(\lambda^n_s+\phi_s)ds \text{ and }\widetilde\mu_1(ds,dx)=\widetilde\mu(dx,ds)-\gamma_s(x)\nu_s(dx)ds.$$

With our Assumptions, we can once more use Girsanov's theorem and define an equivalent probability measure $\P^n$ such that 
$$\frac{d\mathbb P^n}{d\mathbb P}=\mathcal E\left(\int_0^{\cdotp}\left( \lambda^n_s +\phi_s\right)dB_s+\int_0^{\cdotp}\int_E\gamma_s(x)\widetilde\mu(dx,ds)\right)_{t_n}.$$

Taking the conditional expectation on both sides of the above inequality, we obtain
\begin{equation*}
\mathcal E^{g}_{t}\left[X_{t_j} \right] =Y^j_t\leq \mathbb E^{\mathbb P^n}_t\left[X_{t_j}\right] +k(J+1)(t_j-t), \text{  }\P-a.s.,\forall t\in\left[t_{j-1},t_j \right]. 
\end{equation*}

In particular, taking $t=t_{j-1}$ we have
\begin{equation*}
X_{t_{j-1}}\leq \mathcal E^{g}_{t_{j-1}}\left[X_{t_j} \right] \leq \E^{\P^n}_{t_{j-1}}\left[X_{t_j}\right] +k(J+1)(t_j-t_{j-1}), \text{  }\mathbb P-a.s. 
\end{equation*}
Hence $(\widetilde{X}_{t_j})_{j=0..n}$ is a $ \P^n$-submartingale. Define now the following quantities 
$$u_t:=b+k(J+1)t\text{ and } l_t:=a+k(J+1)t.$$

Then, we can apply the classical upcrossing inequality for $\widetilde X$, $u$ and $l$
$$\E^{\P^n}[U_l^u[\widetilde{X},n]]\leq \frac{\E^{\P^n}\left[\left( \widetilde{X}_T-l_T\right)^+ \right]}{u_T-l_T}\leq\frac{\No{X}_{\mathcal{S}^\infty}+2k(J+1)T+\abs{a}}{b-a}.$$

Notice then finally that $U_l^u[\widetilde{X},n]=U_a^b[X,n]$, which implies the desired result.

\vspace{0.5em}
{\bf Step $2$:} When (ii) holds. Using the same arguments as in the proof of (ii) of Proposition \ref{prop.comp}, we can show, using the concavity of $g$ and Assumption \ref{assump:hquad}, that for any $\theta\in(0,1)$
\begin{align*}
\theta g_t(y,z,u)&\leq g_t(0,0,0)+C\theta\abs{y}+(1-\theta)M +\frac{\gamma}{1-\theta}\theta^2\abs{z}^2+\frac{1-\theta}{\gamma}j_t\left(\frac{\gamma\theta u}{1-\theta}\right)\\
&\leq C\theta\abs{y}+(2-\theta)M +\frac{\gamma}{1-\theta}\theta^2\abs{z}^2+\frac{1-\theta}{\gamma}j_t\left(\frac{\gamma\theta u}{1-\theta}\right).
\end{align*}

Hence, considering as in Step $1$ for any $j=0...n$, the solution $Y^j$ of \reff{eq:bsdej_up}, we can use the same exponential transformation as in Step $2$ of the proof of Proposition \ref{prop.comp} to obtain
\begin{align*}
\exp\left(\frac{\gamma\theta}{1-\theta} Y_t^j\right)&\leq \mathbb E_t\left[\exp\left(\frac{\gamma\theta}{1-\theta}X_{t_j}+\gamma \frac{2-\theta}{1-\theta}M(t_j-t)+\frac{\gamma C}{1-\theta}\int_t^{t_j}\abs{Y^j_s}ds\right)\right]\\
&\leq \exp\left(\gamma\left(\frac{C J}{1-\theta}+M(2-\theta)\right)(t_j-t)\right)\mathbb E_t\left[\exp\left(\frac{\gamma\theta}{1-\theta}X_{t_j}\right)\right]\\
&\leq \exp\left(k_\theta(1+J)(t_j-t)\right)\mathbb E_t\left[\exp\left(\frac{\gamma\theta}{1-\theta}X_{t_j}\right)\right],
\end{align*}
for some constant $k_\theta$ depending on $\gamma$, $C$, $M$ and $\theta$.

\vspace{0.5em}
As in Step $1$, choosing $t=t_{j-1}$ and using the fact that $X$ is a $g$-submartingale, we deduce that $(\widehat X_{t_j})_{j=0..n}$ is a $\mathbb P$-submartingale, where
$$\widehat X_t:=\exp\left(k_\theta(1+J)t+\frac{\gamma\theta}{1-\theta}X_t\right).$$

Define now the quantities
$$u^\theta_t:=\exp\left(k_\theta(1+J)t+\frac{\gamma\theta}{1-\theta}b\right)\text{ and }l^\theta_t:=\exp\left(k_\theta(1+J)t+\frac{\gamma\theta}{1-\theta}a\right).$$

We apply the classical upcrossing inequality for $\widehat X$, $u^\theta$ and $l^\theta$
$$ \E[U_{l^\theta}^{u^\theta}[\widehat{X},n]]\leq \frac{\E\left[\left( \widehat{X}_T-l^\theta_T\right)^+ \right]}{u^\theta_T-l^\theta_T}\leq\frac{\exp\left(\frac{\gamma\theta}{1-\theta}\No{X}_{\mathcal{S}^\infty}\right)+\exp\left(\frac{\gamma\theta}{1-\theta}\abs{a}\right)}{\exp\left(\frac{\gamma\theta}{1-\theta}b\right)-\exp\left(\frac{\gamma\theta}{1-\theta}a\right)},$$
which ends the proof, noticing that $U_{l^\theta}^{u^\theta}[\widehat{X},n]=U_{a}^{b}[X,n].$
\ep

\vspace{0.5em}
With this upcrossing inequality in hand, we can argue exactly as in \cite{ma} (see Corollary $5.6$) to obtain
\begin{Corollary}
Let $g$ be as in Theorem \ref{up}. Then any g-sub(super)martingale $X$ admits a c\`adl\`ag modification and furthermore for any countable dense subset $\mathcal D$ of $[0,T]$, it holds for all $t\in[0,T]$ that the two following limits exist $\mathbb P-a.s.$
$$\underset{r\uparrow t,\ r\in\mathcal D}{\lim} X_r \text{ and }\underset{r\downarrow t,\ r\in\mathcal D}{\lim} X_r.$$
\end{Corollary}

\section{Dual Representation and Inf-Convolution}\label{section.app}

We generalize in this section some results of Barrieu and El Karoui \cite{elkarbar2} to the case of quadratic BSDEs with jumps. We give a dual representation of the related $g$-expectations, viewed as convex dynamic risk measures and then we compute in an explicit manner the inf-convolution of two convex $g$-expectations.

\subsection{Dual Representation of the \texorpdfstring{$g$}{g}-expectation}

We will assume in this section that $g_t(y,z,u)=g_t(z,u)$ is independent of $y$ and that the function $g$ is convex. We will prove a dual Legendre-Fenchel type representation for the functional $\Ec^g$, making use of the Legendre-Fenchel transform of $g$. This problem has been treated by Barrieu and El Karoui \cite{elkarbar2} in the case of quadratic BSDEs, we extend it here to the case of quadratic BSDEs with jumps. 

\vspace{0.5em}
In this section, $\Ec^g$ will correspond to a time consistent dynamic convex risk measures. Hence $\Ec^g$ admits a dual representation, as in \cite{elkarbar2}. In this particular case of risk measures constructed from backward SDEs, the penalty function appearing in the dual representation is an integral of the Legendre-Fenchel transform of the generator $g$.
The operator $\Ec^g$, viewed as a time-consistent dynamic convex risk measure has interesting economic applications in insurance. 

\vspace{0.5em}
For $\mu \in \R^d$ and $v \in L^2(\nu_t)$, define the Legendre-Fenchel transform of $g$ in $(z,u)$ as follows
\begin{align*}
G_t(\mu,v):=\sup_{(z,u)\in \R^d\times L^2(\nu_t)} \left\{ \mu.z + \langle v,u\rangle_{t}-g_t(z,u)\right\}.
\end{align*}

Let $\Ac$ denote the space of applications $v \in \mathbb J^2_{\rm{BMO}} \cap L^{\infty}(\nu)$ such that there exists a constant $\delta>0$ with $v_t(x)\geq-1+\delta$, $\P \times dt \times d\nu_t$ -a.e.
\begin{Theorem}\label{thm.rep.duale}
Let $g$ be a given convex function in $(z,u)$ and let Assumptions \ref{assump:hquad} and \ref{assump:hh} hold; assume that $g(0,0)$ and the process $\alpha$ appearing in Assumption \ref{assump:hquad}(iii) are bounded by some constant $M>0$ (then, the existence and uniqueness of the solution of the BSDEJ with generator $g$ and terminal condition $\xi\in L^\infty$ hold by Theorem \ref{th.unique} $(ii)$). We then have
		\begin{enumerate}
			\item[$\rm{(i)}$] For any $\xi_T \in L^{\infty}$,  $$\Ec_t^g(\xi_T) =\underset{(\mu,v) \in \mathbb H^2_{\rm{BMO}} \times \Ac}{\esup^\P} \left\{ \E_t^{\mathbb Q^{\mu,v}}\left[\xi_T- \int_t^T G_s(\mu_s,v_s)ds\right] \right\},\ \mathbb P-a.s.,$$ 
			where $\mathbb Q^{\mu,v}$ is the probability measure defined by $$ \frac{d \mathbb Q^{\mu,v}}{d \P} = \Ec \left( \int_0^. \mu_s dB_s + \int_0^. \int_E v_s(x) \widetilde{\mu}(ds,dx) \right).$$
			\item[$\rm{(ii)}$] Moreover, there exist measurable functions $\overline{\mu}(w,t)$ and $\overline{v}(\omega,t,\cdot)$ such that 
					\begin{align}\Ec_t^g(\xi_T) = \E_t^{\mathbb Q^{\overline{\mu},\overline{v}}}\left[\xi_T- \int_t^T G_s(\overline{\mu}_s,\overline{v}_s)ds\right],\ \mathbb P-a.s. \label{rep_exact1} \end{align}
		\end{enumerate}
\end{Theorem}

\proof
Thanks to the Kazamaki criterion (see for instance Lemma 4.1 in \cite{mor}), we know that if $\mu \in \mathbb H^2_{\rm{BMO}}$ and $v \in \mathbb J^2_{\rm{BMO}}$, then $\Gamma_{\mu,v}:=\frac{d\mathbb  Q^{\mu,v}}{d \P}$ is a true martingale and the probability measure $\mathbb Q^{\mu,v}$ is well defined. $\Ec_t^g(\xi_T)$ is by definition solution of
\begin{align}
\Ec_t^g(\xi_T) =& \ \xi_T + \int_t^T g_s(Z_s,U_s) ds - \int_t^T Z_sdB_s- \int_t^T \int_E U_s(x) \widetilde{\mu}(ds,dx) \nonumber \\
=&\  \xi_T + \int_t^T \left[g_s(Z_s,U_s) -\mu_s.Z_s -\langle v_s,U_s\rangle_{s} \right] ds \nonumber \\
&- \int_t^T Z_sdB^{\mu}_s- \int_t^T \int_E U_s(x) \widetilde{\mu}^{v}(ds,dx), \ \mathbb P-a.s.,\label{inegG}
\end{align}
where $B^{\mu}_t:= B_t - \int_0^t \mu_s ds$ is a $\mathbb Q^{\mu,v}$-Brownian motion and
$$ \widetilde{\mu}^{v}([0,t],A) := \widetilde{\mu}([0,t],A) - \int_0^t \int_A v_s(x) \nu_s(dx) ds \text{  is a  } \mathbb Q^{\mu,v}- \text{martingale}.$$
By Lemma $3.1$ in \cite{kpz_quadratic_1}, $Z \in \mathbb H^2_{\rm{BMO}}$ and $U \in \mathbb J^2_{\rm{BMO}}$. Let us prove that we also have $(Z,J)\in\mathbb H^2(\mathbb Q^{\mu,v})\times\mathbb J^2(\mathbb Q^{\mu,v})$. Indeed, using the number $r>1$ given in Proposition $2.4$ of \cite{kpz_quadratic_1} 
\begin{align*}
\mathbb E^{Q^{\mu,v}}\left[\int_0^T\abs{Z_s}^2ds\right]&=\mathbb E\left[\Gamma_{\mu,v}\int_0^T\abs{Z_s}^2ds\right]\\
&\leq \left(\mathbb E\left[\Gamma_{\mu,v}^r\right]\right)^{1/r}\left(\mathbb E\left[\left(\int_0^T\abs{Z_s}^2ds\right)^q\right]\right)^{1/q}<+\infty,
\end{align*}
where $1/r+1/q=1$ and where we used the energy inequality (($2.3$) in \cite{kpz_quadratic_1}, we refer the reader to \cite{kaz} for more details). The proof for $J$ is the same.
Moreover, 
\begin{align*}
-G_t(\mu,v) = - \underset{(z,u)\in \R^d\times L^2(\nu_t)}{\sup} \left\{ \mu.z + \langle v,u\rangle_{t}-g_t(z,u)\right\} \leq g_t(0,0),
\end{align*}
which means that $-G_t(\mu,v)$ is $\mathbb Q^{\mu,v} \times dt$-integrable. Using these integrability properties and the definition of $G$, we take the conditional expectation in (\ref{inegG}) to obtain 
\begin{align}
\Ec_t^g(\xi_T) &\leq \E_t^{\mathbb Q^{\mu,v}}\left[\xi_T- \int_t^T G_s(\mu_s,v_s)ds \right]. \label{inegG1}
\end{align}
By our assumptions, $g$ is $C^2$ in $z$ and twice Fr\'echet differentiable in $u$, then $\partial g(Z_t,U_t)$ contains a unique element, where the subdifferential $\partial g$ is defined by
\begin{equation*}
\begin{split}
 \partial g(Z_t,U_t) = \Big\{ (\mu,v) \in \R^d \times L^2(\nu_t) \text{, } g_t(z',u') \geq g_t(Z_t,U_t) - \mu.(z'-Z_t) - \langle v,u'-U_t\rangle_{t}, \; \forall (z',u') \Big\}.
\end{split}
\end{equation*}

We take $(\overline{\mu},\overline{v}) \in \partial g(Z_t,U_t)$. We have 
$$g_t(Z_t,U_t) = \overline{\mu}_t.Z_t + \langle \overline{v}_t,U_t\rangle_{t} - G_t(\overline{\mu}_t,\overline{v}_t).$$

We refer to \cite{elkarbar2} for the measurability of $\overline{\mu}$ and $\overline{v}$ with respect to the variable $\omega$. We use Remark \ref{rem_assump_lip} to write
\begin{align*}
 \abs{g_t(z,u)} &\leq \abs{g_t(z,0)} + \No{D_u g_t(z,0)}_{L^2(\nu_t)}^2 \No{u}_{L^2(\nu_t)}^2 + C\No{u}_{L^2(\nu_t)}^2 \\
&\leq \abs{g_t(0,0)} + C\abs{z}^2 + C\left(1+\No{u}_{L^2(\nu_t)}^2\right) \leq C\abs{z}^2 + C\left(1+\No{u}_{L^2(\nu_t)}^2\right),
\end{align*}
where $C$ is a constant whose value may vary from line to line. Putting the above estimation in $G$ leads to
\begin{align*}
 G_t(\overline{\mu}_t,\overline{v}_t) &= \sup_{(z,u)\in \R^d\times L^2(\nu_t)} \left\{\overline{\mu}_t.z + \langle\overline{v}_t,u\rangle_{t} - g_t(z,u) \right\} \\
&\geq  \sup_{u\in L^2(\nu_t)} \left\{ \langle \overline{v}_t,u\rangle_{t} -C -C\No{u}_{L^2(\nu_t)}^2 \right\} + \sup_{z\in \R^d} \left\{\overline{\mu}_t.z -\frac{\gamma}{2}\abs{z}^2 \right\} \\
&= \frac{1}{4C} \No{\overline{v}_t}_{L^2(\nu_t)}^2 -C + \frac{1}{4C} \abs{\overline{\mu}_t}^2.
\end{align*}
From this, we deduce that for $\epsilon <\frac{1}{4C}$,
\begin{align*}
 \left(\frac{1}{4C}-\epsilon\right)\left(\No{\overline{v}_t}_{L^2(\nu_t)}^2 +\abs{\overline{\mu}_t}^2\right) \leq & \ G_t(\overline{\mu}_t,\overline{v}_t) + C -\epsilon \No{\overline{v}_t}_{L^2(\nu_t)}^2 -\epsilon \abs{\overline{\mu}_t}^2 \\
=& \ C -g_t(Z_t,U_t) + \langle\overline{v}_t,U_t\rangle_{t} -\epsilon \No{\overline{v}_t}_{L^2(\nu_t)}^2 +\overline{\mu}_t.Z_t - \epsilon \abs{\overline{\mu}_t}^2 \\
\leq &\ C-g_t(Z_t,U_t) + \frac{1}{4\epsilon}\No{U_t}_{L^2(\nu_t)}^2 + \frac{1}{4\epsilon}\abs{Z_t}^2.
\end{align*}
Since $|g_t(Z_t,U_t)|^{\frac{1}{2}}$ and $U$ are respectively in $\mathbb H^2_{\rm{BMO}}$ and $\mathbb J^2_{\rm{BMO}}$, using the fact that $$\left(\frac{1}{4C}-\epsilon\right)\No{\overline{v}_t}_{L^2(\nu_t)}^2 \leq \left(\frac{1}{4C}-\epsilon\right)\left(\No{\overline{v}_t}_{L^2(\nu)}^2+\abs{\overline{\mu}_t}^2\right) \text{ and,}$$ $$\left(\frac{1}{4C}-\epsilon\right)\abs{\overline{\mu}_t}^2 \leq \left(\frac{1}{4C}-\epsilon\right)\left(\No{\overline{v}_t}_{L^2(\nu_t)}^2+\abs{\overline{\mu}_t}^2\right),$$ 
we obtain that $\overline{v}$ is in $\mathbb J^2_{\rm{BMO}}$ and $\overline{\mu}$ is in $\mathbb H^2_{\rm{BMO}}$. 

\vspace{0.5em}
Furthermore, by our assumptions, $\overline{v} = D_u g \geq -1+ \delta$ and $\overline{v}$ is bounded, then $\overline{v} \in \Ac$. The inequality \reff{inegG1} is thus an equality, and the representation \reff{rep_exact1} holds true.
\ep

\subsection{Inf-Convolution of \texorpdfstring{$g$}{g}-expectations}

Let $g^1_t(z,u)$ and $g^2_t(z,u)$ be two convex functions such that 
\begin{align}
 (g^1 \square g^2)(t,0,0) = \inf_{(\mu,v) \in \R^d \times L^2(\nu_t)} \left\{ g^1_t(\mu,v) + g^2_t(-\mu,-v) \right\} > 0. \label{cond_inf_convol}
\end{align}

The aim of this Section is to compute the optimal risk transfer between two economic agents using $\Ec^{g^1}$ and $\Ec^{g^2}$ as risk measures. The total risk is modeled by a $\Fc_T$-measurable random variable $\xi_T$. The optimal risk transfer will be given through the inf-convolution of the risk measures $\Ec^{g^1}$ and $\Ec^{g^2}$.

\vspace{0.5em}
At time $t$, both agents assess their risk using a monotone convex monetary risk measure (resp. $\Ec_t^{g^1}$ and $\Ec_t^{g^2}$).
For a given loss level $\xi_T$, agent $1$ will take in charge $\xi_T-F$ and transfer to the second agent a quantity $F$, and for this he will pay a premium $\pi(F)$.

\vspace{0.5em}
Agent $1$ minimizes his risk under the constraint that a transaction takes place, he solves : 
\begin{equation}
\inf_{F,\pi} \{ \Ec_t^{g^1}(\xi_T-F + \pi) \} \; \; \text{under the constraint} \; \; \Ec_t^{g^2}(F - \pi) \leq \Ec_t^{g^2}(0)=0  \label{pbm1}
\end{equation}
Binding this last constraint and using the cash-additivity property for $\Ec_t^{g^2}$ gives the optimal price $\pi = \Ec_t^{g^2}(F) - \Ec_t^{g^2}(0) =\Ec_t^{g^2}(F)$.
This is an indifference pricing rule for the first agent, that is to say the price at which he is indifferent (from a risk perspective) between entering and not entering into the transaction. Replacing $\pi = \Ec_t^{g^2}(F)$ in (\ref{pbm1}) and using the cash-additivity property of $\Ec_t^{g^1}$, the insurer program becomes equivalent to the following one:
\[
\inf_{F} \{ \Ec_t^{g^1}(\xi_T-F) + \Ec_t^{g^2}(F)  \} 
=: \Ec_t^{g^1} \square \Ec_t^{g^2} (\xi_T)
\]
We are left with the inf-convolution of $\Ec^{g^1}$ and $\Ec^{g^2}$, problem for which we give some explicit solutions in Theorem \ref{thm.infconvol} and in Section \ref{inf_convol_examples}.

\vspace{0.5em}
More precisely, we will show that, provided that all the quantities considered behave well enough and are in the right spaces, we can identify the inf-convolution of $\Ec^{g^1}$ and $\Ec^{g^2}$ as the solution of a BSDEJ whose generator is the inf-convolution of $g^1$ and $g^2$. Furthermore, we will explicitly construct two $\Fc_T$-measurable random variables $F^{(1)}_T$ and $F^{(2)}_T$ such that $F^{(1)}_T +F^{(2)}_T = \xi_T$ and
\begin{align*}
(\Ec^{g^1} \square \Ec^{g^2})(\xi_T) = \Ec^{g^1}(F^{(1)}_T) + \Ec^{g^2}(F^{(2)}_T).
\end{align*}
We will say that $(F^{(1)}_T, \, F^{(2)}_T)$ is the optimal risk transfer between the agents $1$ and $2$.

\vspace{0.5em}
For this purpose, and for the sake of simplicity, we will assume throughout this section that the solutions to all the considered BSDEJs exist. Notice that this is not such a stringent assumption. Indeed, when it comes to the growth condition of Assumption \ref{assump:hquad}, if we assume that $g^1$ has quadratic growth in $z$ and $u$ and is strongly convex in $(z,u)$, that is to say that there exists some constant $C>0$ such that
$$g^1_t(z,u)- \frac{C}{2}\left(\abs{z}^2+\No{u}_{L^2(\nu_t)}^2\right),$$
is convex, then, since $g^2$ is convex, it is classical that $g^1\square g^2$ also has quadratic growth.

\vspace{0.5em}
Furthermore, we are convinced that as in the classical results by Kobylanski \cite{kob} in the continuous case, this growth condition should be enough to obtain existence of maximal and minimal solutions to the corresponding BSDEJs.

\begin{Remark}
Notice that since the generators are defined on $\Omega\times[0,T]\times\mathbb R^d\times L^2(\nu)\cap L^\infty(\nu)$, a quadratic growth in $u$ is equivalent to the exponential growth assumed in Assumption \ref{assump:hquad}.
\end{Remark}

\begin{Theorem}\label{thm.infconvol}
Let $g^1$ and $g^2$ be two given generators satisfying the assumptions of Theorem \ref{thm.rep.duale}. Denote $(\Ec_t^{1,2}(\xi_T),Z_t,U_t)$ the solution of the BSDEJ with generator $g^1 \square g^2$ and terminal condition $\xi_T$, and let $(\Zh^{(1)}_t,\Uh^{(1)}_t)$ and $(\Zh^{(2)}_t,\Uh^{(2)}_t)$ be four predictable processes such that
\begin{align}
 (g^1 \square g^2)(t,Z_t,U_t) = g^1_t(\Zh^{(1)}_t,\Uh^{(1)}_t) + g^2_t(\Zh^{(2)}_t,\Uh^{(2)}_t) \; dt \times \P-\text{a.s.}
\end{align}
Then, $\rm{(i)}$ For any $\Fc_T$-measurable r.v $F \in L^{\infty}$,
    \begin{align}
     \Ec_t^{1,2}(\xi_T) \leq \Ec_t^{g^1}(\xi_T-F) + \Ec_t^{g^2}(F), \; \P-\text{a.s.}, \; \forall t \in [0,T]. \label{ineg_infconvol}
    \end{align}

$\rm{(ii)}$ Define
 \begin{align*}
     F^{(2)}_T := -\int_0^T g^2_s(\Zh^{(2)}_s,\Uh^{(2)}_s)ds + \int_0^T \Zh^{(2)}_s dB_s + \int_0^T \int_E \Uh^{(2)}_s(x) \widetilde{\mu}(ds,dx),
     \end{align*}
     and assume that the BSDEJs with generators $g^1$ and $g^2$ and terminal conditions $\xi_T-F^{(2)}_T$ and $F^{(2)}_T$ have a solution. If furthermore $\Zh^{(i)} \in \mathbb H^2_{\rm{BMO}}$ and $\Uh^{(i)} \in \mathbb J^2_{\rm{BMO}}$, $i=1,2$, then
    \begin{align}
     (\Ec^{g^1} \square \Ec^{g^2})_t(\xi_T) = \Ec_t^{1,2}(\xi_T) = \Ec_t^{g^1}(\xi_T-F^{(2)}_T) + \Ec_t^{g^2}(F^{(2)}_T). \label{egal_infconvol}
    \end{align}
    \end{Theorem}

\proof
$(\Zh^{(2)},\Uh^{(2)})$ is well defined and predictable thanks to Proposition 8.1 in \cite{elkarbar2}. For $F \in L^{\infty}$, $\Ec_t^{g^1}(\xi_T-F) + \Ec_t^{g^2}(F)$ is solution of
\begin{align*}
 d(\Ec_t^{g^1}(\xi_T-F) + \Ec_t^{g^2}(F)) =& \left( g^1_t(Z^1_t,U^1_t) + g^2_t(Z^2_t,U^2_t)\right)dt -(Z^1_t+Z^2_t)dB_t\\
  &- \int_E (U^1_t(x) + U^2_t(x))\widetilde{\mu}(dt,dx)\\
=&\left( g^1_t(Z_t-Z^2_t,U_t-U^2_t)+g^2_t(Z^2_t,U^2_t)\right)dt -Z_tdB_t - \int_E U_t(x) \widetilde{\mu}(dt,dx)
\end{align*}
and $\Ec_T^{g^1}(\xi_T-F) + \Ec_T^{g^2}(F) = \xi_T$. 

\vspace{0.5em}
Since $g^1$ satisfies the convexity Assumption \ref{assump.bh}, then the function $(z,u) \mapsto g^1_t(z-Z^2_t,u-U^2_t) + g^2_t(Z^2_t,U^2_t)$ is also convex, and we can apply the comparison theorem which directly implies inequality \reff{ineg_infconvol}.

\vspace{0.5em}
Assume now that $\Zh^{(i)} \in \mathbb H^2_{\rm{BMO}}$ and $\Uh^{(i)} \in \mathbb J^2_{\rm{BMO}}$, $i=1,2$, and define $F^{(i)}_t$ by the following forward equations
\begin{align*}
 F^{(1)}_t := \Ec_0^{1,2}(\xi_T) - \int_0^t g^1_s(\Zh^{(1)}_s,\Uh^{(1)}_s)ds + \int_0^t \Zh^{(1)}_s dB_s + \int_0^t \int_E \Uh^{(1)}_s(x) \widetilde{\mu}(ds,dx),\\
 F^{(2)}_t := - \int_0^t g^2_s(\Zh^{(2)}_s,\Uh^{(2)}_s)ds + \int_0^t \Zh^{(2)}_s dB_s + \int_0^t \int_E \Uh^{(2)}_s(x) \widetilde{\mu}(ds,dx).   
\end{align*}
Then we have 
\begin{align*}
 F^{(i)}_t = F^{(i)}_T + \int_t^T g^i_s(\Zh^{(i)}_s,\Uh^{(i)}_s)ds - \int_t^T \Zh^{(i)}_s dB_s - \int_t^T \int_E \Uh^{(i)}_s(x) \widetilde{\mu}(ds,dx),
\end{align*}
and by uniqueness, $F^{(i)}_t = \Ec_t^{g^i}(F^{(i)}_T)$. Moreover, by definition, we have $F^{(1)}_T + F^{(2)}_T = \xi_T$. Since 
\begin{align}
 (g^1 \square g^2)(t,Z_t,U_t) = g^1_t(\Zh^{(1)}_t,\Uh^{(1)}_t) + g^2_t(\Zh^{(2)}_t,\Uh^{(2)}_t) \; dt \times \P-\text{a.s.},
\end{align}
we have the equality 
$$\Ec_t^{1,2}(\xi_T) = \Ec_t^{g^1}(F^{(1)}_T) + \Ec_t^{g^2}(F^{(2)}_T).$$

\vspace{0.5em}
We can conclude that the processes $\Ec_t^{g^1}(F^{(1)}_T) + \Ec_t^{g^2}(F^{(2)}_T)$ and $\Ec_t^{1,2}(\xi_T)$ are solution of the BSDEJ with coefficients $(g^1 \square g^2, \xi_T)$, by uniqueness we have that equality \reff{egal_infconvol} holds.
\ep

\begin{Remark}\label{inf_convol_uniqueness}
The optimal structure $(F^{(1)}_T, F^{(2)}_T)$ is defined up to a constant, more precisely, $(F^{(1)}_T+m, F^{(2)}_T-m)$ with $m\in \R$ is again an optimal structure. Indeed, the cash-additivity property implies that $$\Ec_t^{g^1}(F^{(1)}_T+m) + \Ec_t^{g^2}(F^{(2)}_T-m) = \Ec_t^{g^1}(F^{(1)}_T) + \Ec_t^{g^2}(F^{(2)}_T).$$
\end{Remark}

\subsection{Examples of inf-convolution}\label{inf_convol_examples}
In this Section, we use the previous result on the inf-convolution of $g$-expectations to treat several particular examples.

\subsubsection{Quadratic and Quadratic}
We first study the inf-convolution of two dynamic entropic risk measure. This example is treated by Barrieu and El Karoui \cite{elkarbar2} by a direct method, they find that the optimal risk transfer is \textit{proportional} in the sense that there exists $a \in (0,1)$ such that
\begin{align*}
(\Ec^{g^1} \square \Ec^{g^2})(\xi_T) = \Ec^{g^1}(a \xi_T) + \Ec^{g^2}((1-a)\xi_T).
\end{align*}
We retrieve here this result using Theorem \ref{thm.infconvol}. For this, we first need to study the inf-convolution of the two corresponding generators $g^i$, $i=1,2$
\begin{align}
g^i_t(z,u):=\frac{1}{2\gamma_i}\abs{z}^2+\gamma_i\int_E\left(e^{\frac{u(x)}{\gamma_i}}-1- \frac{u(x)}{\gamma_i}\right)\nu_t(dx),\label{gen.entrop}
\end{align}
where $(\gamma_1,\gamma_2)\in\mathbb R^*_+\times\mathbb R^*_+$.

\vspace{0.5em}
\begin{Lemma}
Let $g^1$ and $g^2$ be the two convex generators defined in equation \reff{gen.entrop}. For any bounded $\Fc_T$-measurable random variable $\xi_T$, we have,
\begin{align*}
(\Ec^{g^1} \square \Ec^{g^2})(\xi_T) = \Ec^{g^1}\left(\frac{\gamma_1}{\gamma_1+\gamma_2} \xi_T\right) + \Ec^{g^2}\left(\frac{\gamma_2}{\gamma_1+\gamma_2} \xi_T\right).
\end{align*}
\end{Lemma}

\proof We can calculate
\begin{align*}
&(g^1\square g^2)(t,z,u)=\underset{v}{\inf}\left\{\frac{1}{2\gamma_1}\abs{v}^2+\frac{1}{2\gamma_2}\abs{z-v}^2\right\}\\
&+\underset{w}{\inf}\left\{\gamma_1\int_E\left(e^{\frac{w(x)}{\gamma_1}}-1- \frac{w(x)}{\gamma_1}\right)\nu_t(dx)+\gamma_2\int_E\left(e^{\frac{u(x)-w(x)}{\gamma_2}}-1- \frac{u(x)-w(x)}{\gamma_2}\right)\nu_t(dx)\right\}.
\end{align*}

The first infimum above is easy to calculate and is attained for $$v^*:=\frac{\gamma_1}{\gamma_1+\gamma_2}z.$$

For the second one, we postulate similarly that it should be attained for
$$w^*:=\frac{\gamma_1}{\gamma_1+\gamma_2}u.$$

In order to verify this result, it is sufficient to prove that for all $(x,y)\in\mathbb R^2$
\begin{equation}\label{gammma}
(\gamma_1+\gamma_2)\left(e^{\frac{x+y}{\gamma_1+\gamma_2}}-1- \frac{x+y}{\gamma_1+\gamma_2}\right)\leq \gamma_1\left(e^{\frac{x}{\gamma_1}}-1- \frac{x}{\gamma_1}\right) + \gamma_2\left(e^{\frac{y}{\gamma_2}}-1- \frac{y}{\gamma_2}\right).
\end{equation}

Set $\lambda:=\gamma_1/(\gamma_1+\gamma_2)$, $a:=x/\gamma_1$, $b:=y/\gamma_2$ and $h(x):=e^x-1-x$, this is equivalent to 
$$h(\lambda a+(1-\lambda)b)\leq \lambda h(a)+(1-\lambda)h(b),$$
which is clear by convexity of the function $h$. Therefore, we finally obtain
$$(g^1\square g^2)(t,z,u)=\frac{1}{2(\gamma_1+\gamma_2)}\abs{z}^2+(\gamma_1+\gamma_2)\int_E\left(e^{\frac{u(x)}{\gamma_1+\gamma_2}}-1- \frac{u(x)}{\gamma_1+\gamma_2}\right)\nu_t(dx).$$

Using the notations of Theorem \ref{thm.infconvol}, we can compute the quantity $F^{(2)}_T$, giving the optimal risk transfer
\begin{align*}
 F^{(2)}_T &= \frac{\gamma_2}{\gamma_1+\gamma_2} \Big(- \int_0^T \left[\frac{1}{2(\gamma_1+\gamma_2)}\abs{Z_t}^2+(\gamma_1+\gamma_2)\int_E\left(e^{\frac{U_t(x)}{\gamma_1+\gamma_2}}-1- \frac{U_t(x)}{\gamma_1+\gamma_2}\right)\nu_t(dx)\right]dt \\
 &+ \int_0^T Z_t dB_t + \int_0^T \int_E U_t(x) \widetilde{\mu}(dt,dx) \Big),\\
&= \frac{\gamma_2}{\gamma_1+\gamma_2} \left(\xi_T-\Ec_0^{1,2}(\xi_T)\right).
\end{align*}

We calculate similarly $F^{(1)}_T$ and obtain
\begin{align*}
 F^{(1)}_T = \frac{\gamma_1}{\gamma_1+\gamma_2}\xi_T + \frac{\gamma_2}{\gamma_1+\gamma_2}\Ec_0^{1,2}(\xi_T).
\end{align*}

\vspace{0.5em}
Now using Remark \ref{inf_convol_uniqueness} with $m=\frac{\gamma_2}{\gamma_1+\gamma_2}\Ec_0^{1,2}(\xi_T)$, we obtain that the proportional structure $\left(\frac{\gamma_1}{\gamma_1+\gamma_2}\xi_T, \frac{\gamma_2}{\gamma_1+\gamma_2}\xi_T\right)$ is optimal.
\ep

\subsubsection{Linear and Quadratic}
Here, we assume that $d=1$. We study the inf-convolution of a dynamic entropic risk measure with a linear one corresponding to a linear BSDEJ. In this case, we want to calculate the inf-convolution of the two corresponding generators $g^1$ and $g^2$ given by
$$g^1_t(z,u):=\frac{1}{2\gamma}\abs{z}^2+\gamma \int_E\left(e^{\frac{u(x)}{\gamma}}-1- \frac{u(x)}{\gamma}\right)\nu_t(dx),$$ \text{ and}
$$g^2_t(z,u):=\alpha z+\beta\int_E (1 \wedge \abs{x})u(x)\nu_t(dx),$$
where $(\gamma,\alpha,\beta)\in\mathbb R^*_+\times\mathbb R\times[-1+\delta,+\infty)$ for some $\delta>0$. 

\begin{Lemma}
Let $g^1$ and $g^2$ be defined in the two previous equations. We have, for any bounded $\Fc_T$-measurable random variable $\xi_T$,
\begin{align*}
(\Ec^{g^1} \square \Ec^{g^2})(\xi_T) = \Ec^{g^1}(F^{(1)}_T) + \Ec^{g^2}(F^{(2)}_T),
\end{align*}
where
\begin{align*}
 F^{(2)}_T = &\xi_T +\frac{1}{2}\alpha^2\gamma T +\gamma\int_0^T \int_E (\beta(1 \wedge \abs{x})-\ln(1+\beta(1 \wedge \abs{x})))\nu_t(dx)dt-\alpha \gamma B_T\\
&- \gamma \int_0^T \int_E \ln(1+\beta(1 \wedge \abs{x}))\widetilde{\mu}(dt,dx),
\end{align*}
and $F^{(1)}_T=\xi_T-F^{(2)}_T$.
\end{Lemma}

\begin{Remark}
Notice that $F^{(2)}_T$ has no longer the linear form with respect to $\xi_T$ obtained in the previous example. Now, the agent 2 receives the value $\xi_T$ perturbed by a random value only depending on the data contained in the filtration, i.e the Brownian motion $B$ and the random measures $\mu$ and $\nu_t$.
\end{Remark}

\proof
We start by computing the inf-convolution in $(z,u)$ of the generators:
\begin{align*}
&(g^1\square g^2)(t,z,u)=\underset{v}{\inf}\left\{\frac{1}{2\gamma}\abs{v}^2+\alpha(z-v)\right\}\\
&+\underset{w}{\inf}\left\{\gamma\int_E\left(e^{\frac{w(x)}{\gamma}}-1- \frac{w(x)}{\gamma}\right)\nu_t(dx)+\beta\int_E (1 \wedge \abs{x})\left(u(x)-w(x)\right)\nu_t(dx)\right\}.
\end{align*}

The first infimum above is easy to calculate and is attained for $$v^*:=\alpha\gamma.$$

Similarly, it is easy to show that the function $w\rightarrow \gamma\left(e^{\frac{w}{\gamma}}-1- \frac{w}{\gamma}\right)+\beta (1 \wedge \abs{x})\left(u(x)-w\right)$ attains its minimum at $w^*:=\gamma\ln(1+\beta(1 \wedge \abs{x})).$ Therefore, we finally obtain
\begin{align*}
(g^1\square g^2)(t,z,u)=&\alpha z-\frac{\alpha^2\gamma}{2}+\int_E\beta(1 \wedge \abs{x})u(x)\nu_t(dx)\\
&+\int_E\gamma\left[\beta(1 \wedge \abs{x})-(1+\beta(1 \wedge \abs{x}))\ln(1+\beta(1 \wedge \abs{x}))\right]\nu_t(dx).
\end{align*}

Notice that all the quantities appearing in $g^1\square g^2$ are finite. Indeed, we first have for any $u\in L^2(\nu)\cap L^\infty(\nu)$
\begin{align}
\abs{\beta(1\wedge\abs{x})u(x)}\leq \beta^2(1\wedge\abs{x}^2)+(u(x))^2,\label{ineg.carre}
\end{align}
and this quantity is therefore $\nu_t$-integrable for all $t$. Then, since $\beta\geq -1+\delta$, it is also clear that for some constant $C_\delta>0$
\begin{align}
0\leq (1+\beta(1 \wedge \abs{x}))\ln(1+\beta(1 \wedge \abs{x}))-\beta(1 \wedge \abs{x})\leq C_\delta(1\wedge\abs{x}^2),\label{ineg.log}
\end{align}
and thus the second integral is also finite.

\vspace{0.5em}
We can now compute for instance the quantity  $F^{(2)}_T$, which is given after some calculations by
\begin{align*}
  F^{(2)}_T =& -\int_0^T \left(\alpha Z_t + \beta \int_E (1 \wedge \abs{x}) U_t(x)\nu_t(dx)\right)dt + \int_0^T Z_t dB_t + \int_0^T\int_E U_t(x)\widetilde{\mu}(dt,dx) \\
&+ \alpha^2 \gamma T +\gamma \beta  \int_0^T \int_E (1 \wedge \abs{x})\ln(1+\beta(1 \wedge \abs{x}))\nu_t(dx)dt -\alpha \gamma B_T\\
& - \gamma \int_0^T \int_E \ln(1+\beta(1 \wedge \abs{x})) \widetilde{\mu}(dt,dx).
\end{align*}
Recall that $(Z,U)$ is part of the solution of the BSDEJ with generator $g^1 \square g^2$ and terminal condition $\xi_T$. Similarly, we can compute the value $ F^{(1)}_T$. Since $ F^{(1)}_T+ F^{(2)}_T =\xi_T$, we obtain
\begin{align*}
 \xi_T = &\frac{1}{2}\alpha^2\gamma T + \gamma\int_0^T \int_E \left((1+\beta(1 \wedge \abs{x}))\ln(1+\beta(1 \wedge \abs{x}))-\beta(1 \wedge \abs{x})\right)\nu_t(dx)dt  \\
&-\int_0^T g^2_t(Z_t,U_t)dt+ \int_0^T Z_t dB_t + \int_0^T\int_E U_t(x)\widetilde{\mu}(dt,dx).
\end{align*}
And finally, we can conclude that the optimal risk transfer takes the form
\begin{align*}
 F^{(2)}_T = &\xi_T +\frac{1}{2}\alpha^2\gamma T +\gamma\int_0^T \int_E (\beta(1 \wedge \abs{x})-\ln(1+\beta(1 \wedge \abs{x})))\nu_t(dx)dt-\alpha \gamma B_T\\
&- \gamma \int_0^T \int_E \ln(1+\beta(1 \wedge \abs{x}))\widetilde{\mu}(dt,dx).
\end{align*}
\ep
\vspace{-0.5em}
\begin{appendix}
\section{Appendix}
\begin{Proposition}\label{concave}
Let $g$ be a function satisfying Assumption \ref{assump:hquad}(i), which is uniformly Lipschitz in $(y,z,u)$. Let $Y$ be a c\`adl\`ag $g$-submartingale (resp. $g$-supermartingale) in $\mathcal S^\infty$. Assume further that $g$ is concave (resp. convex) in $(z,u)$ and that
\begin{align*}
&\hspace{3.4em}g_t(y,z,u)\geq -M-\beta\abs{y}-\frac{\gamma}{2}\abs{z}^2-\frac{1}{\gamma}j_t(-\gamma u).\\
&\left(\text{resp. }g_t(y,z,u)\leq M+\beta\abs{y}+\frac{\gamma}{2}\abs{z}^2+\frac{1}{\gamma}j_t(\gamma u) \right).
\end{align*}
Then there exist a predictable non-decreasing (resp. non-increasing) process $A$ null at $0$ and processes $(Z,U)\in\mathbb H^2\times\mathbb J^2$ such that
$$Y_t=Y_T+\int_t^Tg_s(Y_s,Z_s,U_s)ds-\int_t^TZ_sdB_s-\int_t^T\int_EU_s(x)\widetilde\mu(dx,ds)-A_T+A_t,\ t\in[0,T].$$
\end{Proposition}

\proof
As usual, we can limit ourselves to the $g$-supermartingale case. Let us consider the following reflected BSDEJ with lower obstacle $Y$, generator $g$ and terminal condition $Y_T$
\begin{align}\label{rbsdej}
\nonumber &\widetilde Y_t=Y_T+\int_t^Tg_s(\widetilde Y_s,\widetilde Z_s,\widetilde U_s)ds-\int_t^T\widetilde Z_sdB_s-\int_t^T\int_E\widetilde U_s(x)\widetilde\mu(dx,ds)+A_T-A_t,\ t\in[0,T]\\
\nonumber &\widetilde Y_t\geq Y_t,\ t\in[0,T],\ \ \int_0^T\left(\widetilde Y_{s^-}-Y_{s^-}\right)dA_s=0,
\end{align}
where all the above holds $\mathbb P-a.s.$

\vspace{0.5em}
The existence and uniqueness of a solution consisting of a predictable non-decreasing process $A$ and a triplet $(\widetilde Y,\widetilde Z,\widetilde U)\in\mathcal S^2\times\mathbb H^2\times\mathbb J^2$ follow from the results of \cite{hama2} for instance, since the generator is Lipschitz and the obstacle is c\`adl\`ag and bounded.

\vspace{0.5em}
We will now show that the process $\widetilde Y$ must always be equal to the lower obstacle $Y$, which will provide us the desired decomposition. We proceed by contradiction and assume without loss of generality that $\widetilde Y_0>Y_0$. For any $\eps>0$, we now define the following bounded stopping time
$$\tau_\eps:=\inf\left\{t>0,\ \widetilde Y_t\leq Y_t+\eps,\ \mathbb P-a.s.\right\}\wedge T.$$

By the Skorokhod condition, it is a classical result that the non-decreasing process $A$ never acts before $\tau_\eps$. Therefore, we have for any $t\in[0,\tau_\eps]$
\begin{equation}\label{rbsdej2}
\widetilde Y_t=\widetilde Y_{\tau_\eps}+\int_t^{\tau_\eps}g_s(\widetilde Y_s,\widetilde Z_s,\widetilde U_s)ds-\int_t^{\tau_\eps}\widetilde Z_sdB_s-\int_t^{\tau_\eps}\int_E\widetilde U_s(x)\widetilde\mu(dx,ds),\ \mathbb P-a.s.
\end{equation}

Consider now the BSDEJ on $[0,\tau_\eps]$ with terminal condition $Y_{\tau_\eps}$ and generator $g$ (existence and uniqueness of the solution are consequences of the result of \cite{barles} or \cite{li})
\begin{equation}\label{rbsdej3}
\widehat Y_t=Y_{\tau_\eps}+\int_t^{\tau_\eps}g_s(\widehat Y_s,\widehat Z_s,\widehat U_s)ds-\int_t^{\tau_\eps}\widehat Z_sdB_s-\int_t^{\tau_\eps}\int_E\widehat U_s(x)\widetilde\mu(dx,ds),\ \mathbb P-a.s.
\end{equation}

Notice also that since $Y$ and $g(0,0,0)$ are bounded, $\widehat Y$ and $\widetilde Y$ are also bounded, as a consequence of classical {\it a priori} estimates for Lipschitz BSDEJs and reflected BSDEJs. Then, using the fact that $g$ is convex in $(z,u)$ and that 
$$g_t(y,z,u)\leq M+\beta\abs{y}+\frac{\gamma}{2}\abs{z}^2+\frac{1}{\gamma}j_t(\gamma u),$$
we can proceed as in the Step $2$ of the proof of Proposition \ref{prop.comp} to obtain that for any $\theta\in(0,1)$
\begin{align*}
\delta Y_0&\leq \frac{1-\theta}{\gamma}\ln\left(\mathbb E\left[\exp\left(\gamma\int_0^{\tau_\eps}\left(M+C\abs{\widehat Y_t}+\frac{C}{1-\theta}\abs{\delta Y_s}\right)ds+\frac{\gamma}{1-\theta}\delta Y_{\tau_\eps}\right)\right]\right)\\
&\leq (1-\theta)\ln(C_0)+C\left(\int_0^{\tau_\eps}\No{\delta Y_s}_\infty ds+\No{\delta Y_{\tau_\eps}}_\infty\right),
\end{align*}
where $\delta Y_s:=\widetilde Y_s-\theta\widehat Y_s$ and $C_0$ is some constant which does not depend on $\theta$. Letting $\theta$ go to $1$, and using the fact that by definition $\widetilde Y_{\tau_\eps}-\theta Y_{\tau_\eps}\leq \eps$, we obtain
\begin{equation}\label{rbsdej4}
\widetilde Y_0-\widehat Y_0\leq C\eps +C\int_0^{\tau_\eps}\No{\widetilde Y_s-\widehat Y_s}_\infty ds.
\end{equation}

But since $Y$ is a $g$-supermartingale, we have by definition that $\widehat Y_0\leq Y_0$. Since we assumed that $\widetilde Y_0>Y_0$, this implies that $\widetilde Y_0>\widehat Y_0$. Therefore, we can use Gronwall's lemma in \reff{rbsdej4} to obtain $0<\widetilde Y_0-\widehat Y_0\leq C_1\eps,$ for some constant $C_1>0$, independent of $\eps$. By arbitrariness of $\eps$, this gives us the desired contradiction and ends the proof.
\ep

\end{appendix}

\end{document}